\documentclass[12pt]{article}
\usepackage{geometry,amsmath,amssymb,graphicx,bbm,stmaryrd}
\newcommand{\dueto}[1]{\textup{\textbf{(#1) }}}
\newcommand{\tmem}[1]{{\em #1\/}}
\newcommand{\tmmathbf}[1]{\ensuremath{\boldsymbol{#1}}}
\newcommand{\tmop}[1]{\ensuremath{\operatorname{#1}}}
\newcommand{\tmstrong}[1]{\textbf{#1}}
\newcommand{\tmtextit}[1]{{\itshape{#1}}}
\newcommand{\udots}{{\mathinner{\mskip1mu\raise1pt\vbox{\kern7pt\hbox{.}}\mskip2mu\raise4pt\hbox{.}\mskip2mu\raise7pt\hbox{.}\mskip1mu}}}
\newenvironment{proof}{\noindent\textbf{Proof\ }}{\hspace*{\fill}$\Box$}
\newtheorem{lemma}{Lemma}
\newtheorem{proposition}{Proposition}
\newtheorem{theorem}{Theorem}
\begin{document}

\title{Schur Polynomials\\
and the Yang-Baxter Equation}
\author{Ben Brubaker, Daniel Bump and Solomon Friedberg\\\\
\small Department of Mathematics, MIT, Cambridge MA 02139-4307, USA\\
\small Department of Mathematics, Stanford University, Stanford CA 94305-2125, USA\\
\small Department of Mathematics, Boston College, Chestnut Hill MA 02467-3806, USA}
\maketitle   

\begin{abstract}
  We describe a parametrized Yang-Baxter equation with nonabelian parameter
  group. That is, we show that there is an injective map $g \mapsto R (g)$
  from $\tmop{GL} (2, \mathbbm{C}) \times \tmop{GL} (1, \mathbbm{C})$ to
  $\tmop{End} (V \otimes V)$ where $V$ is a two-dimensional vector space such
  that if $g, h \in G$ then $R_{12} (g) R_{13} (g h) \, R_{23} (h) = R_{23}
  (h) \, R_{13} (g h) R_{12} (g)$. Here $R_{i j}$ denotes $R$ applied to the
  $i, j$ components of $V \otimes V \otimes V$. The image of this map consists
  of matrices whose nonzero coefficients $a_1$, $a_2$, $b_1$, $b_2$, $c_1$,
  $c_2$ are the Boltzmann weights for the six-vertex model, constrained to
  satisfy $a_1 a_2 + b_1 b_2 - c_1 c_2 = 0$. This is the exact center of the
  disordered regime, and is contained within the free Fermionic eight-vertex
  models of Fan and Wu. As an application, we give a new proof based on the
  Yang-Baxter equation of a result of Hamel and King representing a Schur
  polynomial times a deformation of the Weyl denominator as the partition
  function of a six-vertex model. Furthermore, the parameter group can be expanded (within
  the eight-vertex model) to a group having $\tmop{GL} (2) \times \tmop{GL}
  (1)$ as a subgroup of index two. In this expanded context we find a second
  representation of Schur polynomials times a different deformation of the
  Weyl denominator as a partition function. These structures give a
  {\tmem{Yang-Baxter system}} in the sense of Hlavat\'y.
\end{abstract}

Baxter's method of solving lattice models in statistical mechanics is based on
the {\tmem{star-triangle relation}}, which is the identity
\begin{equation}
  \label{shortstar} R_{12} S_{13} T_{23} = T_{23} S_{13} R_{12},
\end{equation}
where $R, S,T$ are endomorphisms of $V \otimes V$ for some
vector space $V$. Here $R_{i j}$ is the endomorphism of $V\otimes V\otimes V$
in which $R$ is applied to the $i$-th and $j$-th
copies of $V$ and the identity map to the $k$-th component, where $i, j, k$
are $1, 2, 3$ in some order. If the endomorphisms $R, S, T$ are all equal, this is
the {\tmem{Yang-Baxter equation}} (cf. {\cite{JimboYB}},
{\cite{MajidQuasitriangular}}).

More generally, one may ask for solutions to a parametrized Yang-Baxter equation,
where the endomorphism $R$ now depends on a parameter $g$ (ranging over a group $G$) 
and (\ref{shortstar}) takes the form
\begin{equation}
  \label{paramyb} R_{12} (g) R_{13} (g \cdot h) R_{23} (h) = R_{23} (h) R_{13} (g \cdot h)
  \, R_{12} (g)
\end{equation}
for arbitrary choice of parameters $g, h \in G$. There are many such examples
in the literature in which the group $G$ is an abelian group such as $\Bbb{R}$
or $\Bbb{R}^\times$. In this paper we present an example of (\ref{paramyb})
having a {\em non-abelian} parameter group. The example arises from
two-dimensional lattice models -- the six- and eight-vertex models.

We now briefly review the connection between lattice models and instances of
(\ref{shortstar}) and (\ref{paramyb}).  In statistical mechanics, one attempts
to understand global behavior of a system from local interactions. To this
end, one defines the partition function of a model to be the sum of certain
locally determined Boltzmann weights over all admissible states of the
system. Baxter (see~{\cite{BaxterInversion}} and~{\cite{Baxter}}, Chapter~9)
recognized that instances of the star-triangle relation allowed one to
explicitly determine the partition function of a lattice model.

The six-vertex, or `ice-type,' model is one such example that is much studied in the literature, and we
revisit it in detail in the next section. For the moment, we offer a few general
remarks needed to describe our results. In our presentation of
the six-vertex model, each state is represented by a labeling of the edges of a finite rectangular
lattice by $\pm$ signs, called {\tmem{spins}}. If the Boltzmann weights are
invariant under sign reversal the system is called {\tmem{field-free}},
corresponding to the physical assumption of the absence of an external field.
For field-free weights, the six-vertex model was solved by Lieb~{\cite{Lieb}}
and Sutherland~{\cite{Sutherland}}, meaning that the partition function
can be exactly computed. The papers of Lieb, Sutherland and Baxter
assume periodic boundary conditions, but non-periodic boundary conditions were
treated by Korepin~{\cite{Korepin}} and Izergin~{\cite{Izergin}}. Much of the
literature assumes that the model is field-free. In this case, Baxter shows 
there is one such parametrized Yang-Baxter 
equation with parameter group $\mathbb{C}^\times$ for each value of a certain 
real invariant $\triangle$, defined below 
in (\ref{baxterdeltadef}) in terms of the Boltzmann weights.


One may ask whether the parameter subgroup $\mathbb{C}^\times$ may be enlarged by including
endomorphisms whose associated Boltzmann weights lie outside the field-free case. 
If $\triangle \neq 0$ the group may {\tmem{not}} be so enlarged. However we will show in
Theorem~\ref{groupcaseagain} that if $\triangle = 0$, then the group $\mathbb{C}^\times$ may
be enlarged to $\tmop{GL} (2, \mathbbm{C}) \times \tmop{GL} (1, \mathbbm{C})$
by expanding the set of endomorphisms to include non-field-free ones. In this
{\tmem{expanded $\triangle = 0$ regime}}, $R (g)$ is not field-free for general
$g$. It is contained within the set of exactly solvable eight-vertex models
called the {\tmem{free Fermionic model}} by Fan and Wu~{\cite{FanWu1}},
{\cite{FanWu2}}. Our calculations suggest that it is not possible to enlarge
the group $G$ to the entire free Fermionic domain in the eight vertex model
but we are able to enlarge $G$ to a group containing $\tmop{GL} (2,
\mathbbm{C}) \times \tmop{GL} (1, \mathbbm{C})$ as a subgroup index two
(Theorem~\ref{groupcaseagainagain}).

In Section~\ref{compo} we give a heuristic argument to show that if there is a
set of endomorphisms such that for any $S$ and $T$ in that set there exists $R$
such that $R_{12} S_{13} T_{23} = T_{23} S_{13} R_{12}$ then an associativity
property is satisfied, so that (\ref{paramyb}) is satisfied. Of course our
rigorous results do not depend on this plausible reasoning, but it seems
useful to know that the associativity that we observe is not entirely
accidental.

As an application of these results, we study the partition function for ice-type models
having boundary conditions determined by an integer partition $\lambda$ and
Boltzmann weights chosen so that both $\triangle = 0$ and so that the degenerate case
$\lambda=0$ matches the standard deformation of Weyl's denominator formula for
$\tmop{GL}_n (\mathbbm{C})$. This leads to an alternate proof of a deformation of the
Weyl character formula for $\tmop{GL}_n$ found by 
Hamel and King~{\cite{HamelKingBijective}}, {\cite{HamelKingUTurn}}. That
result was a substantial generalization of an earlier generating function
identity found by Tokuyama~{\cite{Tokuyama}}, expressed in the language of
Gelfand-Tsetlin patterns.

More precisely, we will exhibit two particular choices of Boltzmann weights and boundary
conditions in the six-vertex model giving systems
$\mathfrak{S}^{\Gamma}_{\lambda}$ and $\mathfrak{S}^{\Delta}_{\lambda}$ for
every partition $\lambda$ of length $\leqslant n$. We will prove that
the partition functions are
\begin{equation}
  \label{thetwopartitionfunctions} Z (\mathfrak{S}_{\lambda}^{\Gamma}) =
  \prod_{i < j} (t_i z_j + z_i) s_{\lambda} (z_1, \cdots, z_n), \hspace{2em} Z
  (\mathfrak{S}_{\lambda}^{\Delta}) = \prod_{i < j} (t_j z_j + z_i)
  s_{\lambda} (z_1, \cdots, z_n),
\end{equation}
where $t_i$ are deformation parameters and $s_{\lambda}$ is the Schur
polynomial (Macdonald~{\cite{Macdonald}}). The method of proof is inspired
by ideas of Baxter in~{\cite{BaxterInversion}} and~{\cite{Baxter}}, though the
Boltzmann weights we use are not field-free. 
The $\Delta$ model is essentially that given by Hamel and
King. The notation here is somewhat unfortunate as $\Delta$ denotes a
recipe for choosing weights and $\triangle$ denotes an invariant defined
in terms of weights, but has been chosen to match earlier uses of this
notation in the literature.

To justify these evaluations of the partition function define
\begin{equation}
  \label{sgammadef} s_{\lambda}^{\Gamma} (z_1, \cdots, z_n ; t_1, \cdots, t_n)
  = \frac{Z (\mathfrak{S}_{\lambda}^{\Gamma})}{\prod_{i < j} (t_i z_j + z_i)}
  .
\end{equation}
Then one seeks to show that $s_{\lambda}^{\Gamma}$ is symmetric in the sense
that it is unchanged if the same permutation is applied to both $z_i$ and
$t_i$. Once this is known, it is possible to show that it is a polynomial in
the $z_i$ and $t_i$, then that it is independent of the $t_i$; finally, taking
$t_i = - 1$ one may invoke the Weyl character formula and conclude that it is
equal to the Schur polynomial.

In order to prove the symmetry property of $s_{\lambda}^{\Gamma}$ we will use
an instance of (\ref{paramyb}) with $\triangle=0$. We thus obtain a new proof of 
Tokuyama's formula and
of Corollary~5.1 in Hamel and King~{\cite{HamelKingBijective}}, which is our
Theorem~\ref{hkthm}. A second instance of the star-triangle relation solves
the same problem for the analogously defined $s_{\lambda}^{\Delta}$, and a
third instance shows directly, without using the above evaluations, that
$s_{\lambda}^{\Gamma} = s_{\lambda}^{\Delta}$.


There are, as we have mentioned, Boltzmann weights of two different types
$\Gamma$ and $\Delta$. (We refer to these as different types of ``ice.'')
Moreover if $X, Y \in \{\Gamma, \Delta\}$ we will give an R-matrix $R_{X Y}$
which has the effect of interchanging a strand of $X$ ice with a strand of $Y$
ice; thus in (\ref{shortstar}), $S$ is of type $X$ and $T$ is of type $Y$. We
will prove that the R-matrices $R_{\Gamma \Gamma}$ and $R_{\Delta \Delta}$
both satisfy the Yang-Baxter equation, and we will prove similar relations
that involve all four types of ice $R_{X Y}$ in various combinations.

Of the six types of ice that we will consider: $\Gamma$, $\Delta$, $R_{\Gamma
\Gamma}$, $R_{\Gamma \Delta}$, $R_{\Delta \Gamma}$ and $R_{\Delta \Delta}$,
only $\Gamma$ and $R_{\Gamma \Gamma}$ come from the space of endomorphisms
parametrized by $\tmop{GL} (2, \mathbbm{C}) \times \tmop{GL} (1,
\mathbbm{C})$. The others may be accommodated by enlarging the parameter group
to a disconnected group having $\tmop{GL} (2, \mathbbm{C}) \times \tmop{GL}
(1, \mathbbm{C})$ as a subgroup of index two.

In another direction, Hlavat\'y~{\cite{HlavatyYBS}} has defined the notion of
a {\tmem{Yang-Baxter system}}. As in our setup, this involves six types of
endomorphisms. His definition has two independent motivations. On the one hand,
there is the work of Freidel and Maillet~{\cite{FreidelMailletQuadratic}} on
integrable systems, and on the other hand, there is work of
Vladimirov~{\cite{VladimirovDoubles}} which attempts to clarify the relation
of the construction of Faddeev, Reshetikhin and Takhtajan~{\cite{FRT}} to
Drinfeld's quantum double. In Section~\ref{ybssection} we show that our
construction is an example of a Yang-Baxter system. In the case where the
$t_i$ are equal, these Yang-Baxter systems are related to those previously
found by Nichita and Parashar~{\cite{NichitaParasharConstructions}},
{\cite{NichitaParasharSpectral}}.

Our boundary conditions depend on the choice of a partition $\lambda$.
Once this choice is made, the states of the model are in bijection with
strict Gelfand-Tsetlin patterns having a fixed top row. These are
triangular arrays of integers with strictly decreasing rows that
interleave (Section~\ref{gtbij}). Since in its original form Tokuyama's
formula expresses what we have denoted $Z (\mathfrak{S}_{\lambda}^{\Gamma})$
as a sum over strict Gelfand-Tsetlin patterns it may be expressed as the
evaluation of a partition function.

This connection between
states of the ice model and strict Gelfand-Tsetlin patterns has one historical
origin in the literature for alternating sign matrices. (An independent
historical origin is in the Bethe Ansatz. See Baxter~\cite{Baxter} Chapter 8
and Kirillov and Reshetikhin~\cite{KirillovReshetikhin}.) The bijection
between the set of alternating sign matrices and strict Gelfand-Tsetlin
patterns having smallest possible top row is in Mills, Robbins and
Rumsey~{\cite{MillsRobbinsRumsey}}, while the connection with what are
recognizably states of the six-vertex model is in Robbins and
Rumsey~{\cite{RobbinsRumsey}}. This connection was used by
Kuperberg~{\cite{Kuperberg}} who gave a second proof (after the purely
combinatorial one by Zeilberger~{\cite{Zeilberger}}) of the alternating
sign matrix conjecture of Mills, Robbins and
Rumsey~{\cite{MillsRobbinsRumsey}}. Kuperberg's paper follows
Korepin~{\cite{Korepin}} and Izergin~{\cite{Izergin}} and makes use of the
Yang-Baxter equation. It was observed by Okada~{\cite{OkadaASM}} and
Stroganov~{\cite{Stroganov}} that the number of $n \times n$ alternating sign
matrices, that is, the value of Kuperberg's ice (with particular Boltzmann
weights involving cube roots of unity) is a special value of the particular
Schur function in $2 n$ variables with $\lambda = (n, n, n - 1, n - 1, \cdots,
1, 1)$ divided by a power of 3. Moreover Stroganov gave a proof using the
Yang-Baxter equation. This occurrence of Schur polynomials in the six-vertex
model is different from the one we discuss, since Baxter's parameter
$\triangle$ is nonzero for these investigations.

There are other works relating symmetric function theory to vertex
models or spin chains. Lascoux~{\cite{Lascoux}}, {\cite{LascouxChern}} gave
six-vertex model representations of Schubert and Grothen\-dieck polynomials of
Lascoux and Sch\"utzenberger~{\cite{LascouxSchutzFlag}} and related these to
the Yang-Baxter equation. Fomin and Kirillov~{\cite{FominKirillov}},
{\cite{FominKirillovGrothendieck}} also gave theories of the Schubert and
Grothendieck polynomials based on the Yang-Baxter equation.
Tsilevich~{\cite{Tsilevich}} gives an interpretation of Schur polynomials and
Hall-Littlewood polynomials in terms of a quantum mechanical system.  Jimbo
and Miwa~\cite{JimboMiwaSolitons}\ give an interpretation of Schur polynomials
in terms of two-dimensional Fermionic systems. (See also
Zinn-Justin~{\cite{ZinnJustin}}.)

McNamara~{\cite{McNamara}} has clarified that the Lascoux papers are
potentially related to ours at least in that the Boltzmann
weights~{\cite{Lascoux}} belong to the expanded $\triangle = 0$ regime.
Moreover, he is able to show based on Lascoux' work how to construct models of
the factorial Schur functions of Biedenharn and Louck.

We are grateful to Gautam Chinta and Tony Licata for stimulating
discussions. This work was supported by NSF grants DMS-0652609, DMS-0652817,
DMS-0652529 and DMS-0702438.  SAGE~{\cite{SAGE}} was very useful in the
preparation of this paper.

\section{The Six-Vertex Model\label{sixver}}\label{sixvertexsection}

We review the six-vertex model from statistical mechanics. Let us consider a
lattice (or sometimes more general graph) in which the edges are labeled with
``spins''~$\pm$. Depending on the spins on its adjacent edges, each vertex
will be assigned a {\tmem{Boltzmann weight}}.

The Boltzmann weight will be zero unless the number of adjacent edges labeled~`$-$' 
is even. Let us denote the possibly nonzero Boltzmann weights as follows:
\[ \begin{array}{|c|c|c|c|c|c|c|c|}
     \hline
     \includegraphics{weighta1.mps} & \includegraphics{weighta2.mps} &
     \includegraphics{weightb1.mps} & \includegraphics{weightb2.mps} &
     \includegraphics{weightc1.mps} & \includegraphics{weightc2.mps} &
     \includegraphics{weightd1.mps} & \includegraphics{weightd2.mps}\\
     \hline
     \includegraphics{rota1.mps} & \includegraphics{rota2.mps} &
     \includegraphics{rotb1.mps} & \includegraphics{rotb2.mps} &
     \includegraphics{rotc1.mps} & \includegraphics{rotc2.mps} &
     \includegraphics{rotd1.mps} & \includegraphics{rotd2.mps}\\
     \hline
     a_1 & a_2 & b_1 & b_2 & c_1 & c_2 & d_1 & d_2\\
     \hline
   \end{array} \]
We will consider the vertices in two possible orientations, as shown above, and
arrange these Boltzmann weights into a matrix as follows:
\begin{equation}
  \label{thisisr} R = \left(\begin{array}{cccc}
    a_1 &  &  & d_1\\
    & b_1 & c_1 & \\
    & c_2 & b_2 & \\
    d_2 &  &  & a_2
  \end{array}\right) = \left(\begin{array}{cccc}
    a_1 (R) &  &  & d_1 (R)\\
    & b_1 (R) & c_1 (R) & \\
    & c_2 (R) & b_2 (R) & \\
    d_2 (R) &  &  & a_2 (R)
  \end{array}\right) .
\end{equation}
If the edge spins are labeled $\nu, \beta, \gamma, \theta \in \{+, -\}$ as
follows:
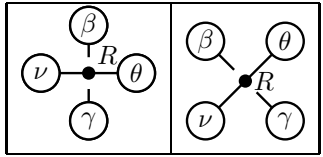
\[ \begin{array}{|c|c|}
     \hline
     \includegraphics{weight0.mps} & \includegraphics{rot0.mps}\\
     \hline
   \end{array} \]
then we will denote by $R_{\nu \beta}^{\theta \gamma}$ the corresponding
Boltzmann weight. Thus $R_{+ +}^{+ +} = a_1 (R)$, etc. 
Because we will sometimes use several different systems of Boltzmann weights
within a single lattice, we label each vertex with the corresponding matrix from which
the weights are taken.

Alternately, $R$ may be thought of as an endomorphism of $V \otimes V$, where $V$
is a two-dimensional vector space with basis $v_+$ and $v_-$. Write
\begin{equation}
  \label{rvinterpretation} R (v_{\nu} \otimes v_{\beta}) = \sum_{\theta,
  \gamma} R_{\nu \beta}^{\theta \gamma} \, v_{\theta} \otimes v_{\gamma} .
\end{equation}
Then the ordering of basis vectors: $v_+ \otimes v_+$, $v_+ \otimes v_-$, $v_-
\otimes v_+$, $v_- \otimes v_-$ gives (\ref{rvinterpretation}) as the matrix~(\ref{thisisr}).

If $\phi$ is an endomorphism of $V \otimes V$ we will denote by $\phi_{12},
\phi_{13}$ and $\phi_{23}$ endomorphisms of $V \otimes V \otimes V$ defined as
follows. If $\phi = \phi' \otimes \phi''$ where $\phi', \phi'' \in \tmop{End}
(V)$ then $\phi_{12} = \phi' \otimes \phi'' \otimes 1$, $\phi_{13} = \phi'
\otimes 1 \otimes \phi''$ and $\phi_{23} = 1 \otimes \phi' \otimes \phi''$. We
extend this definition to all $\phi$ by linearity. Now if $\phi, \psi, \chi$
are three endomorphisms of $V \otimes V$ we define the {\tmem{Yang-Baxter
commutator}}
\[ \left\llbracket \phi, \psi, \chi \right\rrbracket = \phi_{12} \psi_{13}
   \chi_{23} - \chi_{23} \psi_{13} \phi_{12} . \]
\begin{lemma}
  The vanishing of $\left\llbracket R, S, T \right\rrbracket$ is equivalent to
  the star-triangle identity
  \begin{equation}
    \label{startriangleide} \begin{array}{ccccc}
      \sum_{\gamma, \mu, \nu} & 
      \vcenter{\hbox to 120pt{\includegraphics{ybl.mps}}}
      & = & \sum_{\delta,
      \phi, \psi} & 
      \vcenter{\hbox to 120pt{\includegraphics{ybr.mps}}}
    \end{array} \; .
  \end{equation}
  for every fixed combination of spins $\sigma, \tau, \alpha, \beta, \rho,
  \theta$.
\end{lemma}

The term {\tmem{star-triangle identity}} was used by Baxter. The meaning of
equation (\ref{startriangleide}) is as follows.  
 For fixed $\sigma, \tau, \alpha, \beta, \rho, \theta, \mu, \nu, \gamma,$ the
  value or Boltzmann weight of the left-hand side is just the product of the
  Boltzmann weights at the three vertices, that is, $R_{\sigma \tau}^{\nu \mu}
  S_{\nu \beta}^{\theta \gamma} T_{\mu \gamma}^{\rho \alpha}$, and similarly
  the right-hand side. Hence the meaning of (\ref{startriangleide}) is that
  for fixed $\sigma, \tau, \alpha, \beta, \rho, \theta,$
  \begin{equation}
    \label{stimatrix} \sum_{\gamma, \mu, \nu} R_{\sigma \tau}^{\nu \mu} S_{\nu
    \beta}^{\theta \gamma} T_{\mu \gamma}^{\rho \alpha} = \sum_{\delta, \phi,
    \psi} T_{\tau \beta}^{\psi \delta} S_{\sigma \delta}^{\phi \alpha} R_{\phi
    \psi}^{\theta \rho} .
  \end{equation}

\medbreak
\begin{proof}
  Let us apply $\left\llbracket R, S, T \right\rrbracket$ to the vector
  $v_{\sigma} \otimes v_{\tau} \otimes v_{\beta}$. On the one hand by
  (\ref{rvinterpretation})
  \begin{eqnarray*}
    R_{12} S_{13} T_{23} (v_{\sigma} \otimes v_{\tau} \otimes v_{\beta}) & = &
    R_{12} S_{13} \sum_{\psi, \delta} T_{\tau \beta}^{\psi \delta} (v_{\sigma}
    \otimes v_{\psi} \otimes v_{\delta})\\
    & = & R_{12} \sum_{\psi, \delta, \phi, \alpha} S_{\tau \beta}^{\psi
    \delta} T_{\sigma \delta}^{\phi \alpha} (v_{\phi} \otimes v_{\psi} \otimes
    v_{\alpha})\\
    & = & \sum_{\psi, \delta, \phi, \alpha, \theta, \rho} T_{\tau
    \beta}^{\psi \delta} S_{\sigma \delta}^{\phi \alpha} R_{\phi \psi}^{\theta
    \rho} (v_{\theta} \otimes v_{\rho} \otimes v_{\alpha}),
  \end{eqnarray*}
  and similarly
  \[ S_{23} T_{13} R_{12} (v_{\sigma} \otimes v_{\tau} \otimes v_{\beta}) =
     \sum_{\nu, \mu . \theta, \gamma, \rho, \alpha} R_{\sigma \tau}^{\nu \mu}
     S_{\nu \beta}^{\theta \gamma} T_{\mu \gamma}^{\rho \alpha} (_{}
     v_{\theta} \otimes v_{\rho} \otimes v_{\alpha}) . \]
  We see that the vanishing of $\left\llbracket R, S, T \right\rrbracket$ is
  equivalent to (\ref{stimatrix}).
\end{proof}

\medskip

In this section we will be concerned with the {\tmem{six-vertex model}} in
which the weights are chosen so that $d_1 = d_2 = 0$ in the table above. In~{\cite{Baxter}},
Chapter~9, Baxter considered conditions for which, given $S$ and $T$, there exists
a matrix $R$ such that $\left\llbracket R, S, T \right\rrbracket = 0$. We
will slightly generalize his analysis. He considered mainly the
{\tmem{field-free case}} where $a_1 (R) = a_2 (R) = a (R)$, $b_1 (R) = b_2 (R)
= b (R)$ and $c_1 (R) = c_2 (R) = c (R)$. The condition $c_1 (R) = c_2 (R) = c
(R)$ is easily removed, but with no gain in generality. The other two
conditions $a_1 (R) = a_2 (R) = a (R)$, $b_1 (R) = b_2 (R) = b (R)$ are more
serious restrictions.

In the field-free case, let
\begin{equation}
  \label{baxterdeltadef} \triangle (R) = \frac{a (R)^2 + b (R)^2 - c (R)^2}{2
  a (R) \, b (R)}, \hspace{2em} a_1 (R) = a_2 (R) = a (R), \hspace{1em}
  \tmop{etc} .
\end{equation}
Then Baxter showed that given any $S$ and $T$ with $\triangle
(S) = \triangle (T)$, there exists an $R$
such that $\left\llbracket R, S, T \right\rrbracket = 0$.

Generalizing this result to the non-field-free case, we find that there are
not one but two parameters
\begin{eqnarray*}
  \triangle_1 (R) & = & \frac{a_1 (R) a_2 (R) + b_1 (R) b_2 (R) - c_1 (R) c_2
  (R)}{2 a_1 (R) b_1 (R)},\\
  \triangle_2 (R) & = & \frac{a_1 (R) a_2 (R) + b_1 (R) b_2 (R) - c_1 (R) c_2
  (R)}{2 a_2 (R) b_2 (R)} .
\end{eqnarray*}
to be considered.

\begin{theorem}
  \label{groupoidcase}Assume that $a_1 (S)$, $a_2 (S)$, $b_1 (S)$, $b_2 (S)$,
  $c_1 (S)$, $c_2 (S)$, $a_1 (T)$, $a_2 (T)$, $b_1 (T)$, $b_2 (T)$, $c_1 (T)$
  and $c_2 (T)$ are nonzero. Then a necessary and sufficient condition for
  there to exist parameters $a_1 (R)$, $a_2 (R)$, $b_1 (R)$, $b_2 (R)$, $c_1
  (R)$, $c_2 (R)$ such that $\left\llbracket R, S, T \right\rrbracket = 0$
  with $c_1 (R), c_2 (R)$ nonzero is that $\triangle_1 (S) = \triangle_1 (T)$
  and $\triangle_2 (S) = \triangle_2 (T)$.
\end{theorem}

\begin{proof}
  Suppose that $\triangle_1 (S) = \triangle_1 (T)$ and $\triangle_2 (S) =
  \triangle_2 (T)$. Then we may take
  \begin{eqnarray}
    a_1 (R) & = & \frac{b_2 (S) a_1 (T) b_1 (T) - a_1 (S) b_1 (T) b_2 (T) +
    a_1 (S) c_1 (T) c_2 (T)}{a_1 (T)} \nonumber\\
    & = & \frac{a_1 (S) b_1 (S) a_2 (T) - a_1 (S) a_2 (S) b_1 (T) + c_1 (S)
    c_2 (S) b_1 (T)}{b_1 (S)},  \text{\label{a1ieq}}
  \end{eqnarray}
  \begin{eqnarray}
    a_2 (R) & = & \frac{b_1 (S) a_2 (T) b_2 (T) - a_2 (S) b_1 (T) b_2 (T) +
    a_2 (S) c_1 (T) c_2 (T)}{a_2 (T)} \nonumber\\
    & = & \frac{a_2 (S) b_2 (S) a_1 (T) - a_1 (S) a_2 (S) b_2 (T) + c_1 (S)
    c_2 (S) b_2 (T)}{b_2 (S)}  \text{\label{a2ieq}}
  \end{eqnarray}
  
  \begin{equation}
    \text{\label{biieq}} b_1 (R) = b_1 (S) a_2 (T) - a_2 (S) b_1 (T),
    \hspace{2em} b_2 (R) = b_2 (S) a_1 (T) - a_1 (S) b_2 (T),
  \end{equation}
  
  \begin{equation}
    \label{ciieq} c_1 (R) = c_1 (S) c_2 (T), \hspace{2em} c_2 (R) = c_2 (S)
    c_1 (T) .
  \end{equation}
  Using $\triangle_1 (S) = \triangle_1 (T)$ and $\triangle_2 (S) = \triangle_2
  (T)$ it is easy to that the two expressions for $a_1 (R)$ agree, and
  similarly for $a_2 (R)$. One may check that $\left\llbracket R, S, T
  \right\rrbracket = 0$. On the other hand, it may be checked that the
  relations required by $\left\llbracket R, S, T \right\rrbracket = 0$ are
  contradictory unless $\triangle_1 (S) = \triangle_1 (T)$ and
  $\triangle_2 (S) = \triangle_2 (T)$.
\end{proof}

\medskip

In the field-free case, these two relations reduce to a single one, $\triangle
(S) = \triangle (T)$, and it is remarkable that $\triangle (R)$ has the same
value: $\triangle (R) = \triangle (S) = \triangle (T)$.

This equality has important implications for the study of {\tmem{row-transfer matrices}},
one of Baxter's original motivations for introducing the star-triangle relation. Given
Boltzmann weights $a_1 (R), a_2 (R), \cdots$, we associate a $2^n
\times 2^n$ matrix $V (R)$. The entries in this matrix are indexed by pairs
$\alpha = (\alpha_1, \cdots, \alpha_n)$, $\beta = (\beta_1, \cdots, \beta_n)$,
where $\alpha_i, \beta_i \in \{\pm\}$. If $\varepsilon_1, \cdots,
\varepsilon_n \in \{\pm\}$ we may consider the Boltzmann weight of the
configuration:
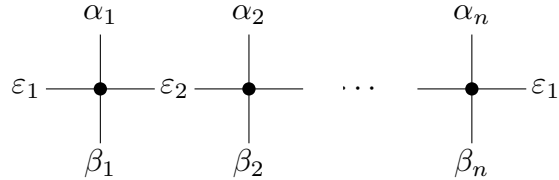
\[ \includegraphics{transfer.mps} \]
Here $\varepsilon_{n + 1} = \varepsilon_1$, so the boundary conditions are
periodic. The coefficient $V (R)_{\alpha, \beta}$ is then the ``partition
function'' for this one-row configuration, that is, the sum over possible
states (assignments of the $\varepsilon_i$).

It follows from Baxter's argument that if $R$ can be found such that
$\left\llbracket R, S, T \right\rrbracket = 0$ then $V (S)$ and $V (T)$
commute, and can be simultaneously diagonalized. We will not review Baxter's
argument here, but variants of it with non-periodic boundary conditions will
appear later in this paper.

In the field-free case when $\left\llbracket R, S, T \right\rrbracket = 0$, $V
(R)$ belongs to the same commuting family as $V (S)$ and $V (T)$. This gives a
great simplification of the analysis in Chapter~9 of Baxter~{\cite{Baxter}}
over the analysis in Chapter~8 using different methods based on the Bethe
Ansatz.

In the non-field-free case, however, the situation is different. If
$\triangle_1 (S) = \triangle_1 (T)$ and $\triangle_2 (S) = \triangle_2 (T)$
then by Theorem~\ref{groupoidcase} there exists $R$ such that $\left\llbracket
R, S, T \right\rrbracket = 0$, and so one may use Baxter's method to prove the
commutativity of $V (S)$ and $V (T)$. However $\triangle_1 (R)$ and
$\triangle_2 (R)$ are not necessarily the same as $\triangle_1 (S) =
\triangle_1 (T)$ and $\triangle_2 (S) = \triangle_2 (T)$, respectively, 
and so $V (R)$ may not commute with $V (S)$ and $V (T)$.

In addition to the field-free case, however, there is {\tmem{another}} case
where $V (R)$ necessarily does commute with $V (S)$ and $V (T)$, and it is
that case which we turn to next. This is the case where $a_1 a_2 + b_1 b_2 -
c_1 c_2 = 0$. The next theorem will show that if the weights of $S$ and $T$
satisfy this condition, then $R$ exists such that
$\left\llbracket R, S, T \right\rrbracket = 0$, and moreover the weights of
$R$ also satisfy the same condition. Thus not only $V(S)$ and $V(T)$ but
also $V(R)$ lie in the same space of commuting transfer matrices.

In this case, with $a_1 = a_1 (R)$, etc., we define
\begin{equation}
  \label{thisispi} \pi (R) = \pi \left(\begin{array}{cccc}
    a_1 &  &  & \\
    & b_1 & c_1 & \\
    & c_2 & b_2 & \\
    &  &  & a_2
  \end{array}\right) = \left(\begin{array}{cccc}
    c_1 &  &  & \\
    & a_1 & b_2 & \\
    & - b_1 & a_2 & \\
    &  &  & c_2
  \end{array}\right) .
\end{equation}
\begin{theorem}
  \label{groupcase}Suppose that
  \begin{equation}
    \label{abcassumption} a_1 (S) a_2 (S) + b_1 (S) b_2 (S) - c_1 (S) c_2 (S)
    = a_1 (T) a_2 (T) + b_1 (T) b_2 (T) - c_1 (T) c_2 (T) = 0.
  \end{equation}
  Then the $R \in End(V \otimes V)$ defined by $\pi (R) = \pi (S) \, \pi (T)^{- 1}$
  satisfies $\left\llbracket R, S, T \right\rrbracket = 0$. Moreover,
  \begin{equation}
    \label{icepreserved} a_1 (R) \, a_2 (R) + b_1 (R) \, b_2 (R) - c_1 (R) \,
    c_2 (R) = 0.
  \end{equation}
\end{theorem}

\begin{proof}
  The matrix $R$ will not be the matrix in Theorem~\ref{groupoidcase}, but
  will rather be a constant multiple of it. We have
  \[ \pi (T)^{- 1} = \frac{1}{D} \left(\begin{array}{cccc}
       c_2 (T) &  &  & \\
       & a_2 (T) & - b_2 (T) & \\
       & b_1 (T) & a_1 (T) & \\
       &  &  & c_1 (T)
     \end{array}\right) \]
  where $D = a_1 (T) a_2 (T) + b_1 (T) b_2 (T) = c_1 (T) c_2 (T)$. With
  notation as in Theorem~\ref{groupoidcase}, using (\ref{abcassumption})
  equations (\ref{a1ieq}) and (\ref{a2ieq}) may be written
  \begin{eqnarray*}
    a_1 (R) & = & a_1 (S) a_2 (T) + b_2 (S) b_1 (T),\\
    a_2 (R) & = & a_2 (S) a_1 (T) + b_1 (S) b_2 (T) .
  \end{eqnarray*}
  Combined with (\ref{biieq}) and (\ref{ciieq}) these imply that $\pi (R) =
  \pi (S) \, D \pi (T)^{- 1}$. However we are free to multiply $R$ by a
  constant without changing the validity of $\left\llbracket R, S, T
  \right\rrbracket = 0$, so we divide it by $D$.
\end{proof}

We started with $S$ and $T$ and produced $R$ such that $\left\llbracket R, S,
T \right\rrbracket = 0$ because this is the
construction motivated by Baxter's method of proving that transfer
matrices commute. However it is perhaps more
elegant to start with $R$ and $T$ and produce $S$ as a function of these. Thus
let $\mathcal{R}$ be the set of endomorphisms $R$ of $V \otimes V$ of the form
(\ref{thisisr}) where $a_1 a_2 + b_1 b_2 = c_1 c_2$. Let $\mathcal{R}^{\ast}$
be the subset consisting of such $R$ such that $c_1 c_2 \neq 0$.

\begin{theorem}
  \label{groupcaseagain}There exists a composition law on $\mathcal{R}^{\ast}$
  such that if $R, T \in \mathcal{R}^{\ast}$, and if $S = R \circ T$ is the
  composition then $\left\llbracket R, S, T \right\rrbracket = 0$. This
  composition law is determined by the condition that $\pi (S) = \pi (R) \pi
  (T)$ where $\pi : \mathcal{R}^{\ast} \longrightarrow \tmop{GL} (4,
  \mathbbm{C})$ is the map (\ref{thisispi}). Then $\mathcal{R}^{\ast}$ is a
  group, isomorphic to $\tmop{GL} (2, \mathbbm{C}) \times \tmop{GL} (1,
  \mathbbm{C})$.
\end{theorem}

\begin{proof}
  This is a formal consequence of Theorem~\ref{groupcase}.
\end{proof}

\medskip

It is interesting that, in the non-field-free case, the group law occurs when $\triangle_1 = \triangle_2 = 0$.
In the application to statistical physics for field-free weights, phase transitions occur when
$\triangle = \pm 1$. If $| \triangle | > 1$ the system is ``frozen'' in the
sense that there are correlations between distant vertices. By contrast $- 1 <
\triangle < 1$ is the disordered range where no such correlations occur, so our group
law occurs in the analog of the middle of the disordered range.

\section{Composition of R-matrices\label{compo}}

Theorem~\ref{groupcaseagain}, defining a group structure on a set of
R-matrices, may be regarded as a non-abelian parametrized Yang-Baxter
equation. In our example, the composition law on R-matrices 
that makes $S$ the product of $R$ and $T$ when $\left\llbracket R, S, T \right\rrbracket = 0$ 
is associative because of its definition in terms of matrix multiplication. In this section, we
give a heuristic argument suggesting that any time we have such a composition
law defined by the vanishing of a Yang-Baxter commutator, associativity should follow. 
This section is not needed for the sequel.

Let us assume that we are given a vector space $V$ over a field $F$ and a
subset $\mathcal{R}$ of $\tmop{End} (V \otimes V)$ which is homogeneous in the
sense that if $0 \neq R \in \mathcal{R}$ then $\mathcal{R}$ contains the
entire ray $F R$. Let $\mathbbm{P}(\mathcal{R})$ be the set of such rays.

Let us assume that if $R, T$ are nonzero elements of $\mathcal{R}$ then there
is another $S \in \mathcal{R}$ that is unique up to scalar multiple such that
$\left\llbracket R, S, T \right\rrbracket = 0$. As we remarked before
Theorem~\ref{groupcase} there might be such an $S$ that would be useable for
applications but that it might not lie in the same space $\mathcal{R}$, and
indeed this is the usual situation for the six vertex model with weights that
are not field-free and also not in the free Fermionic case of
Theorem~\ref{groupcase}. But with this assumption, $(R, T) \mapsto S$ is a
well-defined composition law on $\mathbbm{P}(\mathcal{R})$. Let us denote this
composition $S = R \circ T$. We will give a plausible argument that this
composition law should be associative.

We begin with three nonzero elements $R, S, T$ of $\mathcal{R}$. We will
compare endomorphisms of $V \otimes V \otimes V \otimes V$. In addition to
identities such as $R_{12} (R \circ S)_{13} S_{23} = S_{23} (R \circ S)_{13}
R_{12}$ we will use identities such as $R_{13} T_{24} = T_{24} R_{13}$ which
are true for arbitrary endomorphisms of $V \otimes V$. Let
$$ X_{134} =  (R \circ S)_{13} (R \circ (S \circ T))_{14} T_{34}, \quad X_{134}' = T_{34} (R \circ (S
   \circ T))_{14} (R \circ S)_{13}. $$

First, we have $S_{23} (S \circ T)_{24} X_{134} R_{12}$ equal to
\begin{eqnarray*}
  S_{23} (S \circ T)_{24} (R \circ S)_{13} (R \circ (S \circ T))_{14} T_{34}
  R_{12} & = & \\
  S_{23} (R \circ S)_{13} (S \circ T)_{24} (R \circ (S \circ T))_{14} R_{12}
  T_{34} & = & \\
  S_{23} (R \circ S)_{13} R_{12} (R \circ (S \circ T))_{14} (S \circ T)_{24}
  T_{34} & = & \\
  R_{12} (R \circ S)_{13} S_{23} (R \circ (S \circ T))_{14} (S \circ T)_{24}
  T_{34} & = & \\
  R_{12} (R \circ S)_{13} (R \circ (S \circ T))_{14} S_{23} (S \circ T)_{24}
  T_{34} & = & \\
  R_{12} (R \circ S)_{13} (R \circ (S \circ T))_{14} T_{34} (S \circ T)_{24}
  S_{23} & = & R_{12} X_{134} (S \circ T)_{24} S_{23}. \end{eqnarray*}
Using another string of manipulations, we have $S_{23} (S \circ T)_{24} X_{134}' R_{12}$
equal to
\begin{eqnarray*}
  S_{23} (S \circ T)_{24} T_{34} (R \circ (S \circ T))_{14} (R \circ S)_{13}
  R_{12} & = & \\
  T_{34} (S \circ T)_{24} S_{23} (R \circ (S \circ T))_{14} (R \circ S)_{13}
  R_{12} & = & \\
  T_{34} (S \circ T)_{24} (R \circ (S \circ T))_{14} S_{23} (R \circ S)_{13}
  R_{12} & = & \\
  T_{34} (S \circ T)_{24} (R \circ (S \circ T))_{14} R_{12} (R \circ S)_{13}
  S_{23} & = & \\
  T_{34} R_{12} (R \circ (S \circ T))_{14} (S \circ T)_{24} (R \circ S)_{13}
  S_{23} & = & \\
  R_{12} T_{34} (R \circ (S \circ T))_{14} (R \circ S)_{13} (S \circ T)_{24}
  S_{23} & = & R_{12} X_{134}' (S \circ T)_{24} S_{23}.
\end{eqnarray*}
Now consider an endomorphism $X$ of $V \otimes V \otimes V$ that is
constrained to satisfy
\[ S_{23} (S \circ T)_{24} X_{134} R_{12} = R_{12} X_{134} (S \circ T)_{24}
   S_{23} . \]
This is a linear equation in the matrix coefficients of $X$ in which the
number of conditions exceeds the number of variables. It is reasonable to
assume that if this has a nonzero solution that solution is determined up to
constant multiple. Therefore (up to a constant) we have $X_{134} = X_{134}'$,
that is,
\[ (R \circ S)_{13} (R \circ (S \circ T))_{14} T_{34} = T_{34} (R \circ (S
   \circ T))_{14} (R \circ S)_{13} . \]
Taking the determinant shows that the constant must be a root of unity. If $F
=\mathbbm{R}$ or $\mathbbm{C}$ and $\mathbbm{P}(\mathcal{R})$ is connected,
then by continuity this constant must be 1. This means that $(R \circ (S \circ
T))_{}$ satisfies the definition of $(R \circ S) \circ T$, so at least
plausibly, a composition law defined this way should be associative.

\section{Gamma ice}

Let $z_1, \cdots, z_n$ and $t_1, \cdots, t_n$ be complex numbers with all $z_i
\neq 0$. We will refer to the $z_i$ as {\tmem{spectral parameters}} and the
$t_i$ as {\tmem{deformation parameters}} since these are the roles these
variables will play when we turn to Tokuyama's theorem. Denote
\[ \Gamma (i) = \left(\begin{array}{cccc}
     1 &  &  & \\
     & t_i & z_i (t_i + 1) & \\
     & 1 & z_i & \\
     &  &  & z_i
   \end{array}\right), \hspace{2em} \pi_{\Gamma} (i) =
   \left(\begin{array}{cccc}
     z_i (t_i + 1) &  &  & \\
     & 1 & z_i & \\
     & - t_i & z_i & \\
     &  &  & 1
   \end{array}\right) . \]
Let $\pi_{\Gamma \Gamma} (i, j) = \tmop{const} \times \pi_{\Gamma} (i)
\pi_{\Gamma} (j)^{- 1}$, where it is convenient to take the constant to be
$z_j (t_j + 1)$. It follows from Theorem~\ref{groupcase} that
\begin{equation}
  \label{gagacom} \left\llbracket R_{\Gamma \Gamma} (i, j), \Gamma (i), \Gamma
  (j) \right\rrbracket = 0,
\end{equation}
where $R_{\Gamma \Gamma} (i, j)$ is related to $\pi_{\Gamma \Gamma} (i, j)$ by
the relation (\ref{thisispi}). Concretely,

\begin{equation}
  \label{ggice} R_{\Gamma \Gamma} = \left(\begin{array}{cccc}
    z_j + t_j z_i &  &  & \\
    & t_i z_j - t_j z_i & z_i (t_i + 1) & \\
    & z_j (t_j + 1) & z_i - z_j & \\
    &  &  & z_i + t_i z_j
  \end{array}\right) .
\end{equation}
The six types of vertices corresponding to the non-zero entries of $\Gamma(i)$ and $R_{\Gamma \Gamma}(i,j)$ are given in Table~\ref{tablegam}, together with their Boltzmann weights.

\begin{table}[h]
  
  \[ \begin{array}{|c|c|c|c|c|c|c|}
    \hline
       \begin{array}{c}
         \tmop{Gamma}\\
         \text{Ice}
       \end{array} & \includegraphics{gamma1a.mps} &
       \includegraphics{gamma6a.mps} & \includegraphics{gamma4a.mps} &
       \includegraphics{gamma5a.mps} & \includegraphics{gamma2a.mps} &
       \includegraphics{gamma3a.mps}\\
       \hline
       \text{\begin{tabular}{c}
         Boltzmann\\
         weight
       \end{tabular}} & 1 & z_i & t_i & z_i & z_i (t_i + 1) & 1\\
       \hline
       \text{\begin{tabular}{c}
         Gamma-\\
         Gamma\\
         R-ice
       \end{tabular}} & 
       \vcenter{\hbox to 42pt{\includegraphics{atom1c.mps}}} &
       \vcenter{\hbox to 42pt{\includegraphics{atom6c.mps}}} & 
       \vcenter{\hbox to 42pt{\includegraphics{atom4c.mps}}} &
       \vcenter{\hbox to 42pt{\includegraphics{atom3c.mps}}} & 
       \vcenter{\hbox to 42pt{\includegraphics{atom2c.mps}}} &
       \vcenter{\hbox to 42pt{\includegraphics{atom5c.mps}}}\\
       \hline
       \text{\begin{tabular}{c}
         Boltzmann\\
         weight
       \end{tabular}} & t_j z_i + z_j & t_i z_j + z_i & t_i z_j - t_j z_i &
       z_i - z_j & (t_i + 1) z_i & (t_j + 1) z_j\\
       \hline
     \end{array} \]
  \caption{\label{tablegam}Boltzmann weights for Gamma ice and Gamma-Gamma
  ice.}
\end{table}

\begin{theorem}
  \label{ybesymmetric}The star-triangle identity
 \begin{equation*}
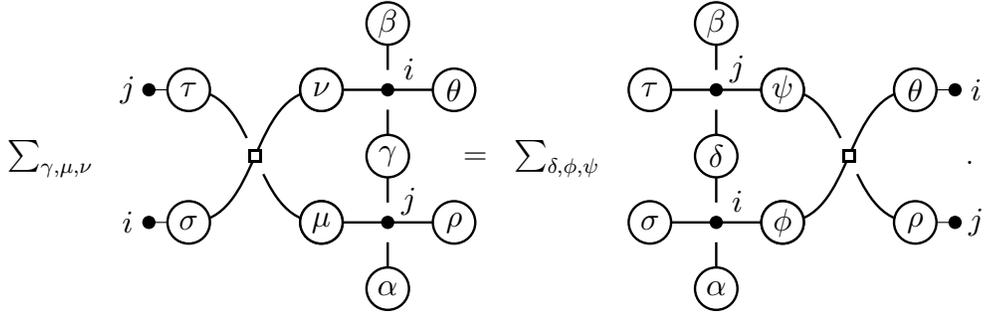
   \begin{array}{ccccc}
      \sum_{\gamma, \mu, \nu} & 
      \vcenter{\hbox to 120pt{\includegraphics{yangbax_g1.mps}}}
      & = & \sum_{\delta,
      \phi, \psi} & 
      \vcenter{\hbox to 120pt{\includegraphics{yangbax_g2.mps}}}
    \end{array} \; .
  \end{equation*}     
is valid with Boltzmann weights as in Table~\ref{tablegam}.
\end{theorem}

\begin{proof}
  This follows from Theorem~\ref{groupcase} since $\pi_{\Gamma \Gamma} (i, j)
  = \tmop{const} \times \pi_{\Gamma} (i) \pi_{\Gamma} (j)^{- 1}$.
\end{proof}

\medskip

We will use Gamma ice to represent Schur polynomials, which are essentially
the characters of finite-dimensional irreducible representations of
$\tmop{GL}_n (\mathbbm{C})$. If $\mu = (\mu_1, \cdots, \mu_n) \in
\mathbbm{Z}^n$ then we may regard $\mu$ as an element of the $\tmop{GL}_n
(\mathbbm{C})$ weight lattice and call it a {\tmem{weight}}. If $\mu_1
\geqslant \cdots \geqslant \mu_n$ we say it is {\tmem{dominant}}, and if \
$\mu_1 > \cdots > \mu_n$ we say it is {\tmem{strictly dominant}}. If $\mu$ is
dominant and $\mu_n \geqslant 0$, it is a {\tmem{partition}}.

\medbreak
{\tmstrong{Note:}} The word ``partition'' occurs in two different senses in
this paper. The partition function in statistical physics is different from
partitions in the combinatorial sense. So for us a reference to a
``partition'' without ``function'' refers to an integer partition. Also
potentially ambiguous is the term ``weight,'' referring to an
element of the $\tmop{GL}_n$ weight lattice, which we identify with
$\mathbbm{Z}^n$. Therefore if we mean Boltzmann weight, we will not omit
``Boltzmann.''

\medbreak
Let $\lambda = (\lambda_1, \cdots, \lambda_n)$ be a fixed partition. We will
denote $\rho = (n - 1, n - 2, \cdots, 0)$. We will consider a rectangular grid
with $n$ rows and $\lambda_1 + n$ columns. We will number the columns of the
lattice in descending order from $\lambda_1 + n - 1$ to $0$.

A {\tmem{state}} of the model will consist of an assignment of ``spins'' $\pm$
to every edge. We will also assign labels to the vertices themselves, which
will be integers between $1$ and $n$. For Gamma ice the vertices in the $i$-th
row will have the label $i$. The spins of the boundary edges are prescribed as
follows.

\medbreak
{\noindent}{\tmstrong{Boundary Conditions determined by $\lambda$.}}
{\tmem{On the left and bottom boundary edges, we put $+$; on the right edges
we put $-$. On the top, we put $-$ at every column labeled $\lambda_i + n - i$
($1 \leqslant i \leqslant n$), that is, for the columns labeled with values in
$\lambda + \rho$. Top edges not labeled by $\lambda_i + n - i$ for any $i$ are
given spin $+$.}}
\medbreak

For example, suppose that $n = 3$ and $\lambda = (3, 1, 0)$, so that $\lambda
+ \rho = (5, 2, 0)$. Then the spins on the boundary are as in the following figure.
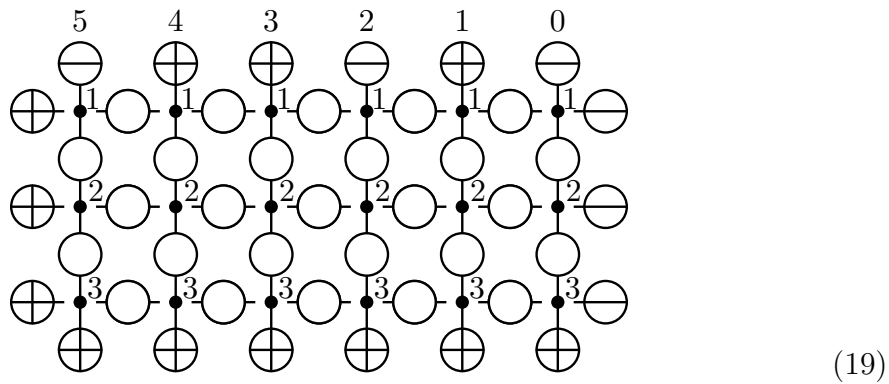
\begin{equation}
  \label{domainwallex} \includegraphics{gamma_ice3.mps}
\end{equation}
The column labels are written at the top, and the vertex labels are written
next to each vertex. The edge spins are marked inside circles. We have left
the edge spins on the interior of the domain blank, since the boundary
conditions only prescribe the spins we have written. The interior spins are
not entirely arbitrary, since we require that at every vertex ``$\bullet$'' the
configuration of spins adjacent to the vertex be one of the six listed in
Table~\ref{tablegam} under ``Gamma ice.''

Let $\mathfrak{S}^{\Gamma}_{\lambda}$ be the {\tmem{Gamma ensemble determined
by $\lambda$}}, by which we mean the set of all such configurations, with the
prescribed boundary conditions. If $x \in \mathfrak{S}^{\Gamma}_{\lambda}$, we
assign a value $w (x)$ called the {\tmem{Boltzmann weight}}. Indeed,
Table~\ref{tablegam} assigns a Boltzmann weight to every vertex, and $w (x)$
is just the product over all the vertices of these Boltzmann weights. The
{\tmem{partition function}} $Z (\mathfrak{S})$ of an ensemble $\mathfrak{S}$
is $\sum_{x \in \mathfrak{S}} w (x)$. As an example, suppose that $n = 2$ and
$l = (0, 0)$ so $\lambda + \rho = (1, 0)$. In this case
$\mathfrak{S}^{\Gamma}_{\lambda}$ has cardinality two, and $Z
(\mathfrak{S}^{\Gamma}_{\lambda}) = t_1 z_2 + z_1$. The states are:
\[ \begin{array}{|c|c|c|}
     \hline
     \text{state} & \includegraphics{state1a.mps} &
     \includegraphics{state2a.mps}\\
     \hline
     \text{Boltzmann weight} & t_1 z_2 & z_1\\
     \hline
   \end{array} \]
The partition function for general $\lambda$ of arbitrary rank will be evaluated 
later in this paper using the star-triangle relation.

%

\section{Gelfand-Tsetlin patterns\label{gtbij}}

Let us momentarily consider a Gamma ice with just one layer of vertices, so
there are three rows of spins. Let $\alpha_1, \cdots, \alpha_m$
be the column numbers (from left to right) of $-$'s in the top row of spins, 
and let $\beta_1, \cdots, \beta_{m'}$ be the column numbers of $-$'s in the
bottom row of edges. For example, in the ice
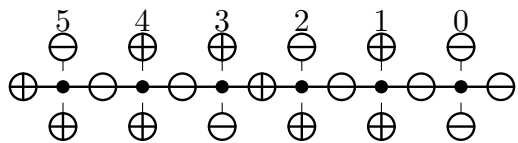
\[ \includegraphics{ice1.mps} \]
we have $m = 3$, $m' = 2$, $(\alpha_1, \alpha_2, \alpha_3) = (5, 2, 0)$ and
$(\beta_1, \beta_2) = (3, 0)$. Since the columns are labeled in decreasing
order, we have $\alpha_1 > \alpha_2 > \cdots$ and $\beta_1 > \beta_2 >
\cdots$.

\begin{lemma}
  \label{onelayerlem}Suppose that the spin at the left edge is $+$. Then we
  have $m = m'$ or $m' + 1$ and $\alpha_1 \geqslant \beta_1 \geqslant \alpha_2
  \geqslant \ldots$ . If $m = m'$ then the spin at the right edge is $+$,
  while if $m = m' + 1$ it is $-$.
\end{lemma}

We express the condition that $\alpha_1 \geqslant \beta_1 \geqslant \alpha_2
\geqslant \ldots$ by saying that the sequences $\alpha_1, \alpha_2, \cdots$
and $\beta_1, \beta_2, \cdots$ {\tmem{interleave}}. This Lemma is essentially
the {\tmem{line-conservation}} principle in Baxter~{\cite{Baxter}}, Section
8.3.

\medbreak
\begin{proof}
  The spins in the middle row are determined by those in the top and bottom
  rows and the left-most spin in the middle row, which is $+$, since the edges
  at each vertex have an even number of $+$ spins. If the rows do not
  interleave then one of the illegal configurations
  \[ \begin{array}{|l|l|}
       \hline
       
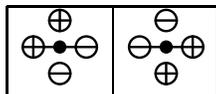
\includegraphics{bad2.mps} & \includegraphics{bad1.mps}\\
       \hline
     \end{array} \]
  will occur. Thus $\alpha_1 \geqslant \beta_1$ since if not, the vertex in
  the $\beta_1$ column would be surrounded by spins in the first illegal
  configuration. Now $\beta_1 \geqslant \alpha_2$ since otherwise the vertex
  in the $\alpha_2$ column would be surrounded by spins in the second above
  illegal configuration, and so forth. The last statement is a consequence of
  the observation that the total number of spins must be even.
\end{proof}

\medskip

We recall that a {\tmem{Gelfand-Tsetlin pattern}} is a triangular array of
dominant weights, in which each row has length one less than the one above it,
and the rows interleave. The pattern is called {\tmem{strict}} if the rows are
strictly dominant.

It follows from Lemma~\ref{onelayerlem} that taking the locations of $-$ in
the rows of vertical lattice edges gives a sequence of strictly dominant weights
forming a strict Gelfand-Tsetlin pattern. For example, given the state
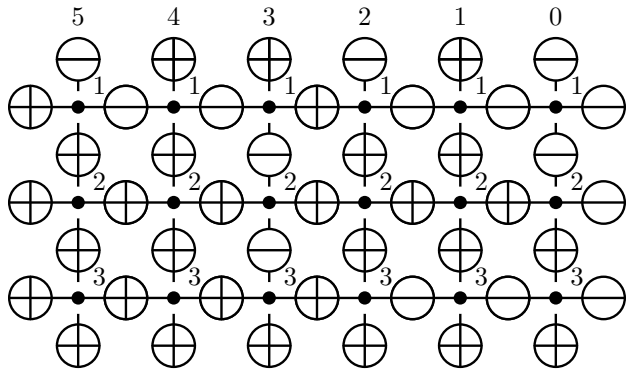
\[ \includegraphics{state4.mps} \]
the corresponding pattern is
\begin{equation}
  \label{examplegtp} \mathfrak{T}= \left\{ \begin{array}{lllll}
    5 &  & 2 &  & 0\\
    & 3 &  & 0 & \\
    &  & 3 &  & 
  \end{array} \right\} .
\end{equation}
It is not hard to see that this gives a bijection between strict
Gelfand-Tsetlin patterns and states with boundary conditions determined
by~$\lambda$. Let us say that the {\tmem{weight}} of a state is $(\mu_1,
\cdots, \mu_n)$ if the Boltzmann weight is the monomial $\tmmathbf{z}^{\mu} =
\prod z_i^{\mu_i}$ times a polynomial in $t_i$. If $\mathfrak{T}$ is a
Gelfand-Tsetlin pattern, let $d_k (\mathfrak{T})$ be the sum of the $k$-th
row. We let $d_{n + 1} (\mathfrak{T}) = 0$.

\begin{lemma}
  \label{weightlem}If $\mathfrak{T}$ is the Gelfand-Tsetlin pattern
  corresponding to a state of weight $\mu$, then $\mu_k = d_k (\mathfrak{T}) -
  d_{k + 1} (\mathfrak{T})$.
\end{lemma}

\begin{proof}
  From Table~\ref{tablegam}, $\mu_k$ is the number of vertices in the $k$-th
  row that have an edge configuration of one of the three forms:
  \[ \begin{array}{|l|l|l|}
       \hline
       \vcenter{\hbox to 42pt
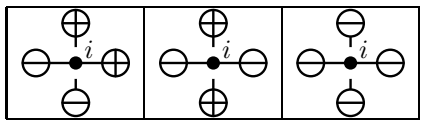
{\includegraphics{gamma2a.mps}}} & 
       \vcenter{\hbox to 42pt{\includegraphics{gamma5a.mps}}} &
       \vcenter{\hbox to 42pt{\includegraphics{gamma6a.mps}}}\\
       \hline
     \end{array} \]
  Let $\alpha_i$'s (respectively $\beta_i$'s) be the column numbers
  for which the top edge spin (respectively, the bottom edge spin) of vertices 
  in the $k$-th row is $-$ (with columns numbered in descending order, as
  always). By Lemma~\ref{onelayerlem} we have $\alpha_1 \geqslant \beta_1
  \geqslant \alpha_2 \geqslant \cdots \geqslant \alpha_{n + 1 - k}$. It is
  easy to see that the vertex in the $j$-column has one of the above
  configurations if and only if its column number $j$ satisfies $\alpha_i > j
  \geqslant \beta_i$ for some $i$. Therefore the number of such $j$ is $\sum
  \alpha_i - \sum \beta_i = d_k (\mathfrak{T}) - d_{k + 1} (\mathfrak{T})$.
\end{proof}

\section{Evaluation of Gamma Ice}

In this section we will prove the following result.

\begin{theorem}
  \label{tokicethm}Let $\lambda = (\lambda_1, \cdots, \lambda_n)$ be a
  partition. Then
  \[ Z (\mathfrak{S}^{\Gamma}_{\lambda}) = \prod_{i < j} (t_i z_j + z_i)
     s_{\lambda} (z_1, \cdots, z_n) . \]
\end{theorem}

To begin with, define
\begin{equation}
  \label{sgamdef} s_{\lambda}^{\Gamma} (z_1, \cdots, z_n ; t_1, \cdots, t_n) =
  \frac{Z (\mathfrak{S}^{\Gamma}_{\lambda})}{\prod_{i < j} (t_i z_j + z_i)} .
\end{equation}
We will eventually show that $s_{\lambda}^{\Gamma}$ is the Schur polynomial
$s_{\lambda}$. But {\tmem{a priori}} it is not obvious from this definition
that $s_{\lambda}^{\Gamma}$ is symmetric, nor that it is a polynomial, nor
that it is independent of~$t$.

\begin{lemma}
  \label{symmetriclabel}The expression $(t_{k + 1} z_k + z_{k + 1}) Z
  (\mathfrak{S}^{\Gamma}_{\lambda})$ is invariant under the interchange of the
  spectral and deformation parameters: $(z_k, t_k) \longleftrightarrow (z_{k + 1}, t_{k +
  1})$.
\end{lemma}

\begin{proof}
  We modify the ice by adding a Gamma-Gamma R-vertex (that is, one of the
  vertices from the bottom row in Table~\ref{tablegam}) to the left of the $k$
  and $k + 1$ rows. Thus (\ref{domainwallex}) becomes (with $k = 2$ for
  illustrative purposes)
  
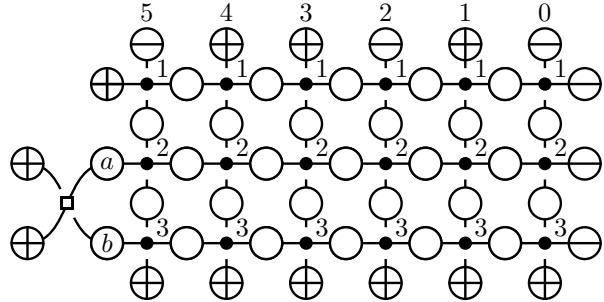
\[ \includegraphics{state5.mps} \]
  which is a new boundary value problem. The only legal values for $a$ and $b$
  are $+$, so every state of this problem determines a unique state of the
  original problem, and the partition function for this state is the original
  partition function multiplied by the Boltzmann weight of the R-vertex, which
  is $t_{k + 1} z_k + z_{k + 1}$. Now we apply the star-triangle identity, and
  obtain equality with the the following configuration.
  
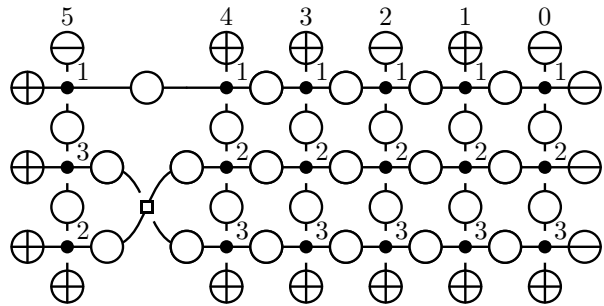
\[ \includegraphics{state6.mps} \]
  Thus if $\mathfrak{S}'$ denotes this ensemble the partition function $Z
  (\mathfrak{S}') = (t_{k + 1} z_k + z_{k + 1}) Z
  (\mathfrak{S}^{\Gamma}_{\lambda})$.
  
  Repeatedly applying the star-triangle identity, we eventually obtain the
  configuration in which the R-vertex is moved entirely to the right.
  
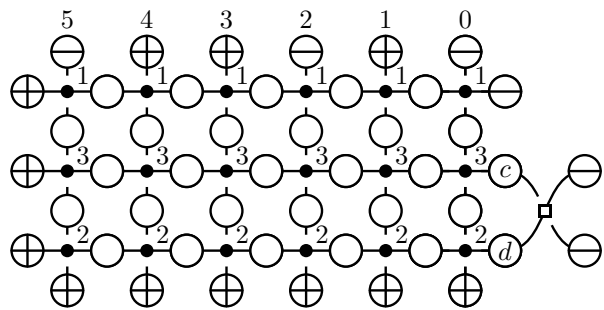
\[ \includegraphics{state7.mps} \]
  Now there is only one legal configuration for the R-vertex, so $c = d = -$.
  The Boltzmann weight at the R-vertex is therefore $t_k z_{k + 1} + z_k$.
  Note that $(z_k, t_k)$ and $(z_{k + 1}, t_{k + 1})$ have been interchanged.
  This proves that $(t_{k + 1} z_k + z_{k + 1}) Z
  (\mathfrak{S}^{\Gamma}_{\lambda})$ is unchanged by switching $(z_k, t_k)$
  and $(z_{k + 1}, t_{k + 1})$.
\end{proof}

\begin{proposition}
  $s_{\lambda}^{\Gamma}$ is a symmetric polynomial in $z_1, \cdots, z_n$, and
  is independent of the~$t_i$.
\end{proposition}

\begin{proof}
  Consider
  \begin{equation}
    \label{partitionnormed} \prod_{i < j} (t_j z_i + z_j) Z
    (\mathfrak{S}^{\Gamma}_{\lambda}) .
  \end{equation}
  We will show that this is invariant under the interchange $k \leftrightarrow
  k + 1$. This means that we interchange both $z_k$ with $z_{k + 1}$ and $t_k$
  with $t_{k + 1}$. Indeed, we may write (\ref{partitionnormed}) as $(t_{k +
  1} z_k + z_{k + 1}) Z (\mathfrak{S}^{\Gamma}_{\lambda})$ times the product
  of all factors $t_j z_i + z_j$ with $i < j$ {\tmem{except}} $(i, j) = (k, k
  + 1)$. These factors are permuted by $k \leftrightarrow k + 1$, so the
  statement follows from Lemma~\ref{symmetriclabel}. Thus
  (\ref{partitionnormed}) is invariant under permutations of the indices,
  where it is understood that the same permutation is applied to the $t_i$ as
  to the $z_i$. Now (\ref{partitionnormed}) equals $\prod_{i \neq j} (t_j z_i
  + z_j) \, s^{\Gamma}_{\lambda} (z_1, \cdots, z_n; t_1, \cdots, t_n),$ so it follows that
  $s^{\Gamma}_{\lambda}$ is also invariant under such permutations. Moreover,
  (\ref{partitionnormed}) is divisible by each $t_j z_i + z_j$ with $i < j$ in
  the unique factorization ring $\mathbbm{C}[z_1, \cdots, z_n, t_1, \cdots,
  t_n]$ . The symmetry property implies that it is also divisible by $t_i z_j
  + z_i$ with $i < j$, and since these are coprime to $\prod_{i < j} (t_j z_i
  + z_j)$ it follows that $Z (\mathfrak{S}^{\Gamma}_{\lambda})$ is divisible
  by these. Therefore $s_{\lambda}^{\Gamma}$ is a polynomial in
  $\mathbbm{C}[z_1, \cdots, z_n, t_1, \cdots, t_n]$.
  
  It remains to be seen that $s_\lambda^\Gamma$ is independent of the $t_i$. In
  \[ s_{\lambda}^{\Gamma} = \frac{Z
     (\mathfrak{S}^{\Gamma}_{\lambda})}{\prod_{i < j} (t_i z_j + z_i)}, \]
  we regard the numerator and the denominator as both being elements of $R
  [t_i]$ where $R =\mathbbm{C}[z_1, \cdots, z_n, t_j (j \neq i)]$. From what
  we have shown, $s^{\Gamma}_{\lambda}$ is a polynomial. We claim that both
  the numerator and denominator have the same degree $i - 1$ in $t_i$. For the
  denominator, this is clear. For the numerator, the number of $-$ in the top
  row of vertical lattice edge spins is $n$ by the boundary conditions, and it follows
  from Lemma~\ref{onelayerlem} that each successive row has one fewer $-$.
  This means that there are $i - 1$ vertices labeled $i$ such that the spin on
  the edge below it is $-$, and from Table~\ref{tablegam}, it follows that the
  number of Boltzmann weights equal to $z_i (t_i + 1)$ or $t_i$ in any
  particular state is $\leqslant i - 1$. The degree of the numerator is thus
  $\leqslant i - 1$ and since the degree of the denominator is $i - 1$, and
  the quotient is a polynomial, both numerator and denominator must have
  degree $i - 1$ in $t_i$. Thus the quotient has degree zero, and does not
  involve $t_i$.
\end{proof}

We may now conclude the proof of Theorem~\ref{tokicethm} by showing that
$s_{\lambda}^{\Gamma} = s_{\lambda}$. Since $s_{\lambda}^{\Gamma}$ is
independent of $t_i$, we may take all $t_i = - 1$. Now in (\ref{sgamdef}) the
denominator becomes $\prod_{i < j} (z_i - z_j)$. Since this is skew-symmetric
under permutations, the numerator $Z (\mathfrak{S}^{\Gamma}_{\lambda})$ is
also skew-symmetric. With $t_i = - 1$ any state containing a vertex in
configuration $\includegraphics{atom8.mps}$ has Boltzmann weight~0, so we are
limited to states omitting this configuration. In view of the bijection
between states and strict Gelfand-Tsetlin patterns, this means that the
corresponding Gelfand-Tsetlin pattern $\mathfrak{T}$ has the property that
every entry from any row but the first is equal to one of the two entries
directly above it. It is easy to see that the weight $\mu$ of such a
coefficient, described by Lemma~\ref{weightlem}, is a permutation $\sigma$ of
the top row of $\mathfrak{T}$, that is, of $\lambda + \rho$. These weights are
all distinct since $\lambda + \rho$ is strongly dominant, i.e. without
repeated entries. Since it is skew-symmetric, its value is $\tmop{sgn}
(\sigma)$ times a constant times $\prod z_j^{\mu_j} = z_j^{\rho_{\sigma (j)} +
\lambda_{\sigma (j)}}$. To determine the constant, we may take the state whose
Gelfand-Tsetlin pattern is
\[ \mathfrak{T}= \left\{ \begin{array}{llllllll}
     \lambda_1 + \rho_1 &  & \lambda_2 + \rho_2 &  &  & \cdots &  &
     \lambda_n\\
     & \lambda_2 + \rho_2 &  &  &  & \lambda_n &  & \\
     &  & \ddots &  & \udots &  &  & \\
     &  &  & \lambda_n &  &  &  & 
   \end{array} \right\} . \]
This has weight $\prod z_j^{\lambda_j + \rho_j}$ and so
\[ s_{\lambda}^{\Gamma} (z_1, \cdots, z_n) = \frac{\sum_{\sigma \in S_n}
   \tmop{sgn} (\sigma) \prod z_j^{\rho_{\sigma (j)} + \lambda_{\sigma
   (j)}}}{\prod_{i < j} (z_j - z_i)} \]
which equals $s_{\lambda} (z_1, \cdots, z_n)$ by the Weyl character formula.

\section{Tokuyama's theorem}

We recall some definitions from Tokuyama {\cite{Tokuyama}}. An entry of a
Gelfand-Tsetlin pattern (not in the top row) is classified as
{\tmem{left-leaning}} if it equals the entry above it and to the left. It is
{\tmem{right-leaning}} if it equals the entry above it and to the right. It is
{\tmem{special}} if it is neither left- nor right-leaning. Thus in
(\ref{examplegtp}), the 3 in the bottom row is left-leaning, the 0 in the
second row is right-leaning and the 3 in the middle row is special. If
$\mathfrak{T}$ is a Gelfand-Tsetlin pattern, let $l (\mathfrak{T})$ be the
number of left-leaning entries. Let $d_k (\mathfrak{T})$ be the sum of the
$k$-th row of $\mathfrak{T}$, and $d_{n + 1} (\mathfrak{T}) = 0$.

\begin{theorem}
  {\dueto{Tokuyama}}We have
  \[ \sum_{\mathfrak{T}} \left( \prod_{k = 1}^n z_k^{d_k (\mathfrak{T}) - d_{k
     + 1} (\mathfrak{T})} \right) t^{l (\mathfrak{T})} (t + 1)^{s
     (\mathfrak{T})} = \prod_{i < j} (z_i + t z_j) s_{\lambda} (z_1, \cdots,
     z_n), \]
  where the sum is over all strict Gelfand-Tsetlin patterns with top row
  $\lambda + \rho$.
\end{theorem}

\begin{proof}
  If $\mathfrak{T}$ corresponds to a state of the Gamma ice with boundary
  conditions determined by $\lambda$, then we will show that the Boltzmann
  weight of the state is the term on the left-hand side. From
  Lemma~\ref{weightlem} the powers of $z$ are correct. It is easy to see that
  if an entry in the $k$-th row of $\mathfrak{T}$ is left leaning (respectively
  special), and that entry is $j$, then the configuration in the $j$-column and the
  $k$-th row of the ice is
  \[ \text{$
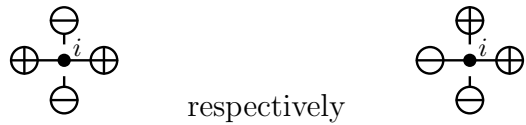
\includegraphics{gamma4a.mps} \hspace{2em} \tmop{respectively}
     \hspace{2em} \includegraphics{gamma2a.mps}$} \]
  so from Table~\ref{tablegam}, it follows that the powers of $t_i$ are also
  correct. The statement now follows from Theorem~\ref{tokicethm}.
\end{proof}

\section{The Yang-Baxter equation for Gamma-Gamma ice\label{gammagamma}}

We will prove a star-triangle relation that only involves Gamma-Gamma ice. Let
us think of Gamma ice as being organized into strands of horizontal lattice edges,
with every Gamma vertex of the strand having the same label $i$. We may think
of Gamma-Gamma ice as a tool that switches two strands. The following result
states that this tool respects the braid relation. We have drawn this picture
differently from that in Theorem~\ref{ybesymmetric} since this Yang-Baxter
equation involves only horizontal edges, while that in
Theorem~\ref{ybesymmetric} involves both horizontal and vertical edges.

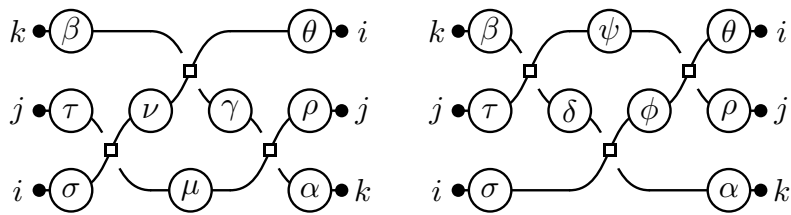
\[ \begin{array}{lll}
     \includegraphics{rice1.mps} &  & \includegraphics{rice2.mps}
   \end{array} \]
With $\sigma, \tau, \beta, \alpha, \rho, \theta$ fixed, we may regard these
two configurations as ensembles each involving three Gamma-Gamma vertices. The
Yang-Baxter equation says that they have the same partition function.

\begin{theorem}
  \label{gammagammaybe}The Yang-Baxter equation is true in the form
  \[ \sum_{\mu, \nu, \gamma} R (j, k)_{\mu \gamma}^{\rho \alpha} R (i, k)_{\nu
     \beta}^{\theta \gamma} R (i, j)_{\sigma \tau}^{\nu \mu} = \sum_{\delta,
     \phi, \psi} R (j, k)_{\tau \beta}^{\psi \delta} R (i, k)_{\sigma
     \delta}^{\phi \alpha} R (i, j)_{\phi \psi}^{\theta \rho}, \]
  with $R = R_{\Gamma \Gamma}$.
\end{theorem}

\medbreak
\begin{proof}
  This follows from Theorem~\ref{groupcaseagain} since $\pi_{\Gamma \Gamma}
  (i, j) = \tmop{const} \times \pi_{\Gamma} (i) \pi_{\Gamma} (j)^{- 1}$, so
  \[ \pi_{\Gamma \Gamma} (i, j) \pi_{\Gamma \Gamma} (j, k) = \tmop{const}
     \times \pi_{\Gamma \Gamma} (i, k) . \]
\end{proof}

\section{More Star-Triangle Relations}

There are further star-triangle relations which go outside the six-vertex
model. We find that the discussion in Section~\ref{sixvertexsection} can be
extended the set of Boltzmann weights in the eight vertex model that has
either $a_1 a_2 + b_1 b_2 - c_1 c_2 = 0$ and $d_1 = d_2 = 0$ or $a_1 a_2 + b_1
b_2 - d_1 d_2 = 0$ and $c_1 = c_2 = 0$. The parameter subgroup will have the
$\tmop{GL} (2, \mathbbm{C}) \times \tmop{GL} (1, \mathbbm{C})$ of
Theorem~\ref{groupcaseagain} as a subgroup of index two. Let
${\widehat{\mathcal{R}^{}}^{\ast}}$ be the set of $R$ as in (\ref{thisisr})
with such weights, where it is assumed $a_1 a_2 + b_1 b_2 \neq 0$.

\begin{theorem}
  \label{groupcaseagainagain}There exists a composition law on
  $\hat{\mathcal{R}}^{\ast}$ such that if $R, T \in \hat{\mathcal{R}}^{\ast}$,
  and if $S = R \circ T$ is the composition then $\left\llbracket R, S, T
  \right\rrbracket = 0$. This composition law is determined by the condition
  that $\pi (S) = \pi (R) \pi (T)$ where $\pi : \hat{\mathcal{R}}^{\ast}
  \longrightarrow \tmop{GL} (4, \mathbbm{C})$ is the map defined by
  (\ref{thisispi}) if $c_1, c_2$ are nonzero, and by
  \[ \pi (R) = \pi \left(\begin{array}{cccc}
       a_1 &  &  & d_1\\
       & b_1 &  & \\
       &  & b_2 & \\
       d_2 &  &  & a_2
     \end{array}\right) = \left(\begin{array}{cccc}
       &  &  & d_1\\
       & i a_2 & - i b_1 & \\
       & i b_2 & i a_1 & \\
       d_2 &  &  & 
     \end{array}\right) \]
  if $d_1, d_2$ are nonzero.
\end{theorem}

Here $i = \sqrt{- 1}$.

\begin{proof}
  Let us call $R \in \hat{\mathcal{R}}^{\ast}$ of {\tmem{Type $C$}} if $c_1,
  c_2$ are nonzero (so $d_1 = d_2 = 0$) and of {\tmem{Type D}} in the other
  case. There are four cases to consider. One, where $R$ and $T$ are both of
  type $C$, is already in Theorem~\ref{groupcaseagain}. In the other three
  cases, we compute $\left\llbracket R, S, T \right\rrbracket = 0$ with $S$ as
  follows.
  
  If $R$ is of type $C$ and $T$ is of type $D$ then $S$ is of type $D$ with
  \begin{eqnarray*}
    a_1 (S) & = & a_2 (R) a_1 (T) + b_1 (R) b_1 (T),\\
    a_2 (S) & = & a_1 (R) a_2 (T) + b_2 (R) b_2 (T),\\
    b_1 (S) & = & - b_2 (R) a_1 (T) + a_1 (R) b_1 (T),\\
    b_2 (S) & = & - b_1 (R) a_2 (T) + a_2 (R) b_2 (T),\\
    d_1 (S) & = & c_1 (R) d_1 (T),\\
    d_2 (S) & = & c_2 (R) d_1 (T) .
  \end{eqnarray*}
  If $R$ is of type $D$ and $T$ is of type $C$ then $S$ is of type $D$ with
  \begin{eqnarray*}
    a_1 (S) & = & a_1 (R) a_2 (T) + b_2 (R) b_2 (T),\\
    a_2 (S) & = & a_2 (R) a_1 (T) + b_1 (R) b_1 (T),\\
    b_1 (S) & = & b_1 (R) a_1 (T) - a_2 (R) b_2 (T),\\
    b_2 (S) & = & b_2 (R) a_1 (T) - a_1 (R) b_1 (T),\\
    d_1 (S) & = & d_1 (R) c_2 (T),\\
    d_2 (S) & = & d_2 (R) c_1 (T) .
  \end{eqnarray*}
  Finally, if $R$ and $T$ are of type $D$ then $S$ is of type $C$ with
  \begin{eqnarray*}
    a_1 (S) & = & - a_2 (R) a_2 (T) + b_1 (R) b_2 (T),\\
    a_2 (S) & = & - a_1 (R) a_1 (T) - b_2 (R) b_1 (T),\\
    b_1 (S) & = & b_2 (R) a_2 (T) + a_1 (R) b_2 (T),\\
    b_2 (S) & = & b_1 (R) a_1 (T) + a_2 (R) b_1 (T),\\
    c_1 (S) & = & d_1 (R) d_2 (T),\\
    c_2 (S) & = & d_2 (R) d_1 (T) .
  \end{eqnarray*}
  These computations may be translated into the identity $\pi (S) = \pi (R) \,
  \pi (T)$.
\end{proof}

We will give some applications of this. The Boltzmann weights for a variety of
other models are given in Table~\ref{tablerest}. While Gamma ice is of Type C
in the terminology of the last proof, we also introduce Delta ice which is of
Type D. Delta-Delta ice is of Type D and Gamma-Delta and Delta-Gamma ice are
of Type C. We will distinguish between Gamma ice and Delta ice by using
$\bullet$ to represent Gamma ice and $\circ$ to represent Delta ice, and
variants of this convention will also distinguish the other four types of ice.

\begin{table}[h]
  
  \[ \begin{array}{|l|l|l|l|l|l|l|}
       \hline
       \text{$\begin{array}{c}
         \tmop{Delta}\\
         \text{Ice}
       \end{array}$} & 
       \vcenter{\hbox to 42pt{\includegraphics{delta1a.mps}}} &
       \vcenter{\hbox to 42pt{\includegraphics{delta2a.mps}}} & 
       \vcenter{\hbox to 42pt{\includegraphics{delta3a.mps}}} &
       \vcenter{\hbox to 42pt{\includegraphics{delta4a.mps}}} & 
       \vcenter{\hbox to 42pt{\includegraphics{delta5a.mps}}} &
       \vcenter{\hbox to 42pt{\includegraphics{delta6a.mps}}}\\
       \hline
       \text{\begin{tabular}{c}
         Boltzmann\\
         weight
       \end{tabular}} & z_i & z_i (t_i + 1) & 1 & z_i t_i & 1 & 1\\
       \hline
       \text{\begin{tabular}{c}
         Delta-\\
         Delta\\
         R-ice
       \end{tabular}} & 
       \vcenter{\hbox to 42pt{\includegraphics{atom1b.mps}}} &
       \vcenter{\hbox to 42pt{\includegraphics{atom2b.mps}}} & 
       \vcenter{\hbox to 42pt{\includegraphics{atom3b.mps}}} &
       \vcenter{\hbox to 42pt{\includegraphics{atom4b.mps}}} & 
       \vcenter{\hbox to 42pt{\includegraphics{atom5b.mps}}} &
       \vcenter{\hbox to 42pt{\includegraphics{atom6b.mps}}}\\
       \hline
       \text{\begin{tabular}{c}
         Boltzmann\\
         weight
       \end{tabular}} & t_i z_i + z_j & z_j (t_j + 1) & t_j z_j - t_i z_i &
       z_i - z_j & (t_i + 1) z_i & z_i + t_j z_j\\
       \hline
       \text{\begin{tabular}{c}
         Gamma-\\
         Delta\\
         R-ice
       \end{tabular}} & 
       \vcenter{\hbox to 42pt{\includegraphics{atom1a.mps}}} &
       \vcenter{\hbox to 42pt{\includegraphics{atom2a.mps}}} & 
       \vcenter{\hbox to 42pt{\includegraphics{atom3a.mps}}} &
       \vcenter{\hbox to 42pt{\includegraphics{atom4a.mps}}} & 
       \vcenter{\hbox to 42pt{\includegraphics{atom5a.mps}}} &
       \vcenter{\hbox to 42pt{\includegraphics{atom6a.mps}}}\\
       \hline
       \text{\begin{tabular}{c}
         Boltzmann\\
         weight
       \end{tabular}} & t_i t_j z_j - z_i & (t_j + 1) z_j & t_i z_j + z_i &
       t_j z_j + z_i & (t_i + 1) z_i & z_i - z_j\\
       \hline
       \text{\begin{tabular}{c}
         Delta-\\
         Gamma\\
         R-ice
       \end{tabular}} & 
       \vcenter{\hbox to 42pt{\includegraphics{atom1d.mps}}} &
       \vcenter{\hbox to 42pt{\includegraphics{atom2d.mps}}} & 
       \vcenter{\hbox to 42pt{\includegraphics{atom3d.mps}}} &
       \vcenter{\hbox to 42pt{\includegraphics{atom4d.mps}}} & 
       \vcenter{\hbox to 42pt{\includegraphics{atom5d.mps}}} &
       \vcenter{\hbox to 42pt{\includegraphics{atom6d.mps}}}\\
       \hline
       \text{\begin{tabular}{c}
         Boltzmann\\
         weight
       \end{tabular}} & z_i - z_j & (t_i + 1) z_i & t_j z_i + z_j & t_i z_i +
       z_j & (t_j + 1) z_j & - t_i t_j z_i + z_j\\
       \hline
     \end{array} \]
  \caption{\label{tablerest}Boltzmann weights for various types of ice with
  spectral parameters $(z_i, t_i)$ and $(z_j, t_j)$. (See Table~\ref{tablegam}
  for Gamma and Gamma-Gamma ice.)}
\end{table}

Thus in addition to (\ref{ggice}) we have:

\[ R_{\Delta \Delta} (z_i, t_i, z_j, t_j) = \left(\begin{array}{cccc}
     z_i t_i + z_j &  &  & \\
     & z_i - z_j & z_j t_j + z_j & \\
     & z_i t_i + z_i & z_j t_j - z_i t_i & \\
     &  &  & z_j t_j + z_i
   \end{array}\right), \]

\[ R_{\Gamma \Delta} (z_i, t_i, z_j, t_j) = \left(\begin{array}{cccc}
     - z_i + t_i t_j z_j &  &  & z_i t_i + z_i\\
     & z_j t_j + z_i &  & \\
     &  & z_j t_i + z_i & \\
     z_j t_j + z_j &  &  & z_i - z_j
   \end{array}\right), \]
\[ R_{\Delta \Gamma} (z_i, t_i, z_j, t_j) = \left(\begin{array}{cccc}
     z_i - z_j &  &  & z_j t_j + z_j\\
     & z_i t_i + z_j &  & \\
     &  & z_i t_j + z_j & \\
     z_i t_i + z_i &  &  & z_j - t_i t_j z_i
   \end{array}\right) . \]
We will denote by $\Gamma (z_i, t_i)$ what was previously denoted $\Gamma
(i)$. We have also $\Delta (z_i, t_i)$:
\[ \Gamma (z_i, t_i) = \left(\begin{array}{cccc}
     1 &  &  & \\
     & t_i & (t_i + 1) z_i & \\
     & 1 & z_i & \\
     &  &  & z_i
   \end{array}\right), \hspace{2em} \Delta (z_i, t_i) =
   \left(\begin{array}{cccc}
     z_i &  &  & 1\\
     & z_i t_i &  & \\
     &  & 1 & \\
     z_i (t_i + 1) &  &  & 1
   \end{array}\right) . \]
\begin{theorem}
  \label{gencirstar}If $X, Y \in \{\Gamma, \Delta\}$ then
  \begin{equation}
    \label{genrxyrel} \left\llbracket R_{X Y} (z_i, t_i, z_j, t_j) \,, X (z_i,
    t_i), Y (z_j, t_j) \right\rrbracket = 0.
  \end{equation}
\end{theorem}

\begin{proof}
  In each of the four cases
  \[ \pi (R_{X Y} (z_i, t_i, z_j, t_j)) \pi (Y (z_j, t_j)) = z_j (t_j + 1)
     \times \pi (X (z_i, t_i)) . \]
  The result then follows from Theorem~\ref{groupcaseagainagain}.
\end{proof}

Now we turn to generalizations of the Yang-Baxter equation. For every choice
of $z$ and $t$ and $X \in \{\Gamma, \Delta\}$, let $V^X (z, t)$ be a
two-dimensional vector space with basis $v_+^X (z, t)$ and $v_-^Y (z, t)$.
Then $R^{X Y} (z_1, t_1, z_2, t_2)$ defines an endomorphism of $V^X (z_1, t_1)
\otimes V^Y (z_2, t_2)$ by
\[ R (v_{\sigma} \otimes v_{\tau}) = \sum_{\mu,\nu} R_{\sigma \tau}^{\nu \mu}
   v_{\nu} \otimes v_{\mu}, \hspace{2em} R = R^{X Y} (z_1, t_1, z_2, t_2) . \]
\begin{theorem}
  \label{parybethm}If $X, Y, Z \in \{\Gamma, \Delta\}$ then we have
  \[ \left\llbracket R_{X Y} (z_1, t_1, z_2, t_2), R_{X Z} (z_1, t_1, z_3,
     t_3), R_{Y Z} (z_2, t_2, z_3, t_3) \right\rrbracket = 0. \]
  Moreover
  \[ \left\llbracket R_{X Y} (z_2, t_1, z_1, t_2), R_{X Z} (z_3, t_1, z_1,
     t_3), R_{Y Z} (z_3, t_2, z_2, t_3) \right\rrbracket = 0. \]
\end{theorem}

\begin{proof}
  This follows from Theorem~\ref{groupcaseagainagain}.
\end{proof}

We now describe the boundary conditions for Delta ice in the ensemble
$\mathfrak{S}_{\lambda}^{\Delta}$ that appears in the second identity in
(\ref{thetwopartitionfunctions}). The columns are labeled, as with the Gamma
ice, in decreasing order. However we label the vertices in decreasing row
order, so the labels of the vertices of the top row are $n$, and so forth.

The Delta ice boundary conditions are as follows. We again fix a partition
$\lambda$. On the left boundary edges, we put $-$; on the right and bottom
edges we put $+$. On the top, we put $-$ at every column labeled $\lambda_i +
n - i$ ($1 \leqslant i \leqslant n$), that is, for the columns labeled with
values in $\lambda + \rho$. Top edges not labeled by $\lambda_i + n - i$ for
any $i$ are given spin $+$. Thus if $\lambda = (3, 1, 0)$, here is the Delta
ice. (To indicate that this is Delta ice, the vertices are marked $\circ$.)

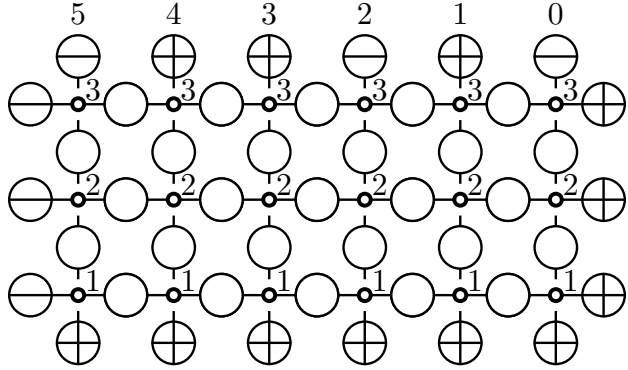
\[ \includegraphics{delta_ice1b.mps} \]
\begin{theorem}
  \label{hkthm}The partition function is
  \[ Z (\mathfrak{S}^{\Delta}_{\lambda}) (z_1, \cdots, z_n ; t_1, \cdots, t_n)
     = \prod_{i < j} (t_j z_j + z_i) s_{\lambda} (z_1, \cdots, z_n) . \]
  
\end{theorem}

\begin{proof}
  This is proved analogously to Theorem~\ref{tokicethm}, using the case $X = Y
  = \Delta$ of Theorem \ref{gencirstar}. We leave the details of the proof to
  the reader.
\end{proof}

Theorem \ref{gencirstar} may be used to show that
\begin{equation}
  \label{statementb} \prod_{i < j} (t_j z_j + z_i) Z
  (\mathfrak{S}^{\Gamma}_{\lambda}) = Z (\mathfrak{S}^{\Delta}_{\lambda})
  \prod_{i < j} (t_i z_j + z_i)
\end{equation}
directly without invoking Theorems~\ref{tokicethm} and~\ref{hkthm}. This fact
is closely related to Statement~B in Brubaker, Bump and
Friedberg~{\cite{wmd5book}}, and the following argument may be used to give an
alternative proof of that result in the special case where the degree
(denoted~$n$ in~{\cite{wmd5book}}) equals~1.

Begin with an element $x$ of $\mathfrak{S}^{\Gamma}_{\lambda}$, say (for
example with $\lambda = (3, 1, 0)$):
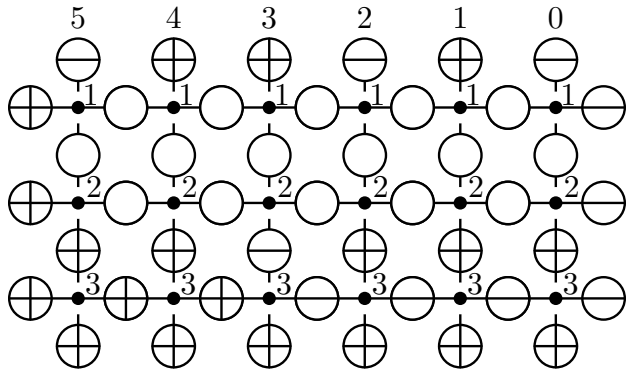
\[ \includegraphics{gamma_ice2a.mps} \]
(The unlabeled edges can be filled in arbitrarily.) We wish to transform this
into an element of an ensemble that has a row of Delta ice so that we may use
the mixed star-triangle relation. We simply change the signs of all the
entries on the edges in the 3 row:
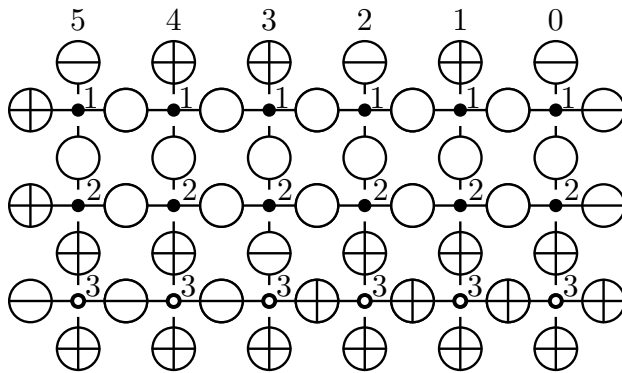
\[ \includegraphics{mixed_ice1a.mps} \]
Let $x'$ be this element of the mixed ensemble $\mathfrak{S}'$. We observe
that the Boltzmann weights satisfy $w (x) = w (x')$. Indeed, in the bottom row
only the following types of Gamma ice can appear:
\[ \begin{array}{|l|l|l|l|}
     \hline
     \begin{array}{c}
       \tmop{Gamma}\\
       \text{Ice}
     \end{array} & \includegraphics{gamma1a.mps} &
     \includegraphics{gamma3a.mps} & \includegraphics{gamma5a.mps}\\
     \hline
     & 1 & 1 & z_i\\
     \hline
   \end{array} \]
These change to:
\[ \begin{array}{|l|l|l|l|}
     \hline
     \text{$\begin{array}{c}
       \tmop{Delta}\\
       \text{Ice}
     \end{array}$} & \includegraphics{gamma5a.mps} &
     \includegraphics{delta3a.mps} & \includegraphics{delta1a.mps}\\
     \hline
     & 1 & 1 & z_i\\
     \hline
   \end{array} \]

Observe that the weights are unchanged. {\tmem{Note that this would not work
in any row but the last}} because it is essential that there be no $-$ on the
bottom edge spins. Now we add a Gamma-Delta R-vertex.
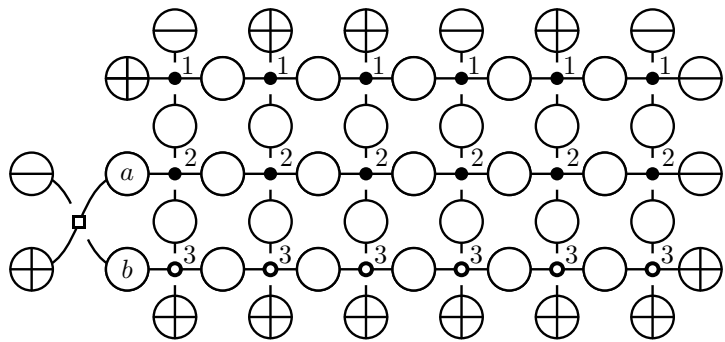
\[ \includegraphics{braided_ice1a.mps} \]
If $\mathfrak{S}''$ is this ensemble, we claim that $Z (\mathfrak{S}'') = (t_3
z_3 + z_2) Z (\mathfrak{S}') = (t_3 z_3 + z_2) Z
(\mathfrak{S}^{\Gamma}_{\lambda})$. Indeed, from Table~\ref{tablerest}, the
values of $a$ and $b$ must be $+, -$ respectively and so the value of the
R-vertex is $t_3 z_3 + z_2$ for every element of the ensemble. Now using the
star-triangle relation, we obtain $Z (\mathfrak{S}'') = Z (\mathfrak{S}''')$
where $\mathfrak{S}'''$ is the ensemble:
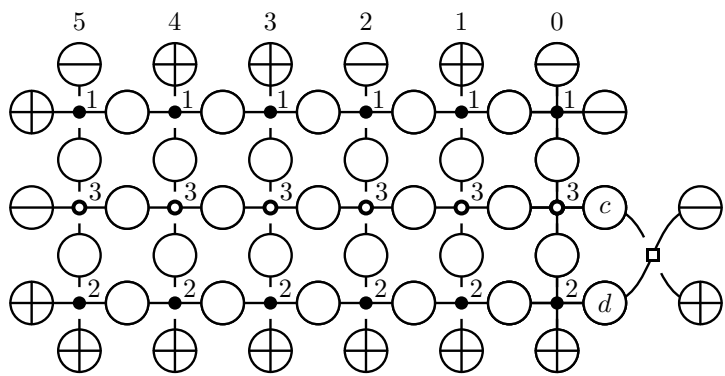
\[ \includegraphics{braided_ice3a.mps} \]
Here we must have $c, d = +, -$and so $(t_3 z_3 + z_2) Z
(\mathfrak{S}^{\Gamma}_{\lambda}) = Z (\mathfrak{S}''') = (t_2 z_3 + z_2) Z
(\mathfrak{S}^{(\tmop{iv})})$ where $\mathfrak{S}^{(\tmop{iv})}$ is the
ensemble:
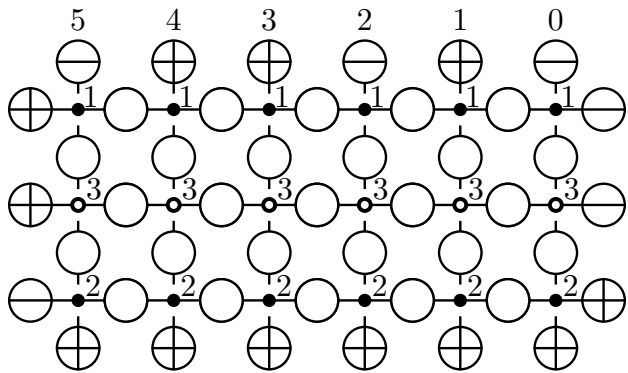
\[ \includegraphics{mixed_ice2a.mps} \]
We repeat the process, first moving the Delta layer up to the top, then
introducing another Delta layer at the bottom, etc., until we have the
ensemble $\mathfrak{S}^{\Delta}_{\lambda}$, obtaining~(\ref{statementb}).

\section{Yang-Baxter Systems}\label{ybssection}

The results of this section are further applications of
Theorem~\ref{groupcaseagainagain}.

An important property of the R-matrices $R_{X Y} (z_i, t_i, z_j, t_j)$ is that
they are {\tmem{projectively triangular}}. That is,
\begin{equation}
  \label{triangular} R_{X Y} (z_i, t_i, z_j, t_j)^{- 1} = c_{X Y} (z_i, t_i,
  z_j, t_j) P \, R_{Y X} (z_j, t_j, z_i, t_i) \, P
\end{equation}
where $c_{X Y} (z_i, t_i, z_j, t_j)$ is a scalar and
\[ P = \left(\begin{array}{cccc}
     1 &  &  & \\
     &  & 1 & \\
     & 1 &  & \\
     &  &  & 1
   \end{array}\right) . \]
The constant $c_{X Y}$ may be eliminated by multiplying $R_{X Y}$ by a suitable
scalar - for example in the case $X = Y = \Gamma$ if $R'_{\Gamma \Gamma} (z_i,
t_i, z_j, t_j) = (z_j t_i + z_i)^{- 1} R_{\Gamma \Gamma} (z_i, t_i, z_j, t_j)$
then $R'_{\Gamma \Gamma}$ satisfies (\ref{triangular}) without the $c_{X Y}$, at the cost of
introducing denominators.

{\tmem{Yang-Baxter systems}} occur with varying degrees of generality in
connection with different problems. One type occurs in the work of
Vladimirov~{\cite{VladimirovDoubles}} on quantum doubles; another type occurs
in Hlavat\'y~{\cite{HlavatyQuantizedBraidedGroups}} on quantized braided groups.
The most general formulation~{\cite{HlavatyNonultralocal}},
{\cite{HlavatyYBS}} involves four types of matrices which correspond to our
$R_{X Y}$, $X, Y \in \{\Gamma, \Delta\}$.

The axioms for a parametrized (or ``colored'') Yang-Baxter system in the most
general definition require four types of matrices, $A, B, C, D$, depending on
parameters $z_1$ and $z_2$ and subject to the properties
\begin{equation}
  \label{ybsaxioms} \begin{array}{ccccccc}
    \left\llbracket A, A, A \right\rrbracket & = & 0, & \hspace{2em} &
    \left\llbracket D, D, D \right\rrbracket & = & 0,\\
    \left\llbracket A, C, C \right\rrbracket & = & 0, &  & \left\llbracket D,
    B, B \right\rrbracket & = & 0,\\
    \left\llbracket A, B^{\ddagger}, B^{\ddagger} \right\rrbracket & = & 0, &
    \hspace{2em} & \left\llbracket D, C^{\ddagger}, C^{\ddagger}
    \right\rrbracket & = & 0,\\
    \left\llbracket A, C, B^{\ddagger} \right\rrbracket & = & 0, &  &
    \left\llbracket D, B, C^{\ddagger} \right\rrbracket & = & 0,
  \end{array}
\end{equation}
where we now denote
\[ \left\llbracket X, Y, Z \right\rrbracket = X_{12} (z_1, z_2) Y_{13} (z_1,
   z_3) Z_{23} (z_2, z_3) - Z_{23} (z_2, z_3) Y_{13} (z_1, z_3) X_{12} (z_1,
   z_2) \]
and $X^{\ddagger} (z_1, z_2) = P X (z_2, z_1) P$. We have two spectral
parameters $z$ and $t$, so we interpret
\[ X^{\ddagger} (z_1, t_1, z_2, t_2) = P X (z_2, t_2, z_1, t_1) P. \]
\begin{theorem}
  Let $X, Y \in \{\Gamma, \Delta\}$. Then
  \[ A = R_{X X}, \hspace{2em} C = B^{\ddagger} = R_{X Y}, \hspace{2em} D =
     {R_{Y Y}}^{\ddagger} \]
  is a Yang-Baxter system satisfying~(\ref{ybsaxioms}).
\end{theorem}

\begin{proof}
  We leave the verification to the reader.
\end{proof}

Note that by projective triangularity we may replace $B$ by ${R_{Y X}}^{- 1}$,
which is a scalar multiple of ${R_{X Y}}^{\ddagger}$. Thus if $X = \Gamma, Y =
\Delta$ we have the Yang-Baxter system
\[ A = R_{\Gamma \Gamma}, \hspace{2em} B = {R_{\Delta \Gamma}}^{- 1},
   \hspace{2em} C = R_{\Gamma \Delta}, \hspace{2em} D = {R_{\Delta
   \Delta}}^{\ddagger}, \]
which uses each of the four braided ice types in Table~\ref{tablerest} exactly once.
It is probably most interesting to take $X \neq Y$, but worth noting that we
can also make a Yang-Baxter system with $R_{\Gamma \Gamma}$ (or
$R_{\Delta \Delta}$) playing all four roles. And we also obtain a Yang-Baxter
system as follows by interchanging the $z_i$ (but not the $t_i$) in the
spectral parameters.

\begin{theorem}
  Another set of four Yang-Baxter systems may be obtained by taking
  \[ A = \hat{R}_{X X}, \hspace{2em} C = B^{\ddagger} = \hat{R}_{X Y},
     \hspace{2em} D = \hat{R}_{Y Y}^{\ddagger}, \]
  where
  \[ \hat{R}_{X Y} (z_1, t_1, z_2, t_2) = R_{X Y} (z_2, t_1, z_1, t_2) . \]
\end{theorem}

\begin{proof}
  We leave this to the reader.
\end{proof}

\end{document}